\documentclass{article}

\usepackage{microtype}
\usepackage{graphicx}
\usepackage{subfigure}
\usepackage{booktabs} 

\usepackage{hyperref}

\usepackage[accepted]{icml2019blank}

\usepackage{graphicx} 
\usepackage{subfigure} 
\usepackage{multirow}
\usepackage{mathtools}
\usepackage{relsize}
\usepackage{url}
\usepackage{nicefrac}
\usepackage{xspace}
\usepackage{algorithm}
\usepackage{algorithmic}
\usepackage{multicol}
\usepackage{float}

\usepackage[utf8]{inputenc} 
\usepackage[T1]{fontenc}    
\usepackage{hyperref}       
\usepackage{url}            

\usepackage{nicefrac}       
\usepackage{microtype}      
 
\usepackage{layout}
\usepackage{multirow}

\usepackage{amsfonts}
\usepackage{amsmath}
\usepackage{amsthm}
\usepackage{amssymb} 

\usepackage{mdframed} 
\usepackage{thmtools}

\usepackage[font=normalsize, labelfont=bf]{caption}

\usepackage{color}
\usepackage{graphicx}

\DeclareMathOperator{\Ocal}{\mathcal{O}} 

\newcommand{\R}{\mathbb{R}} 

\newcommand{\cO}{{\cal O}}

\newcommand{\mA}{{\bf A}}

\newcommand{\mP}{{\bf P}}

\usepackage[colorinlistoftodos,bordercolor=orange,backgroundcolor=orange!20,linecolor=orange,textsize=scriptsize]{todonotes}

\newcommand{\Sam}{S}

\newcommand{\eqdef}{:=}
\newcommand{\norm}[1]{\lVert#1\rVert}      

\DeclareMathOperator{\Prob}{Prob}

\newcommand{\Diag}[1]{\mathbf{Diag}\left( #1\right)}

\providecommand{\trace}[1]{{\rm Trace}\left( #1\right)}

\newcommand{\Exp}[1]{{{\rm E}}\left[#1\right] }    
\newcommand{\E}[1]{{\rm E}\left[#1\right] }

\theoremstyle{plain}
\newtheorem{theorem}{Theorem}  
\newtheorem{lemma}[theorem]{Lemma} 

\theoremstyle{definition}
\newtheorem{definition}{Definition}

\theoremstyle{remark}
\icmltitlerunning{Nonconvex Variance Reduced Optimization with Arbitrary  Sampling} 

\begin{document}

\twocolumn[
\icmltitle{Nonconvex Variance Reduced  Optimization  with Arbitrary  Sampling\textsuperscript{*}}



\icmlsetsymbol{equal}{*}

\begin{icmlauthorlist}
\icmlauthor{Samuel Horv\'{a}th}{kaust}
\icmlauthor{Peter Richt\'{a}rik}{kaust,mipt,ed}
\end{icmlauthorlist}

\icmlaffiliation{kaust}{King Abdullah University of Science and Technology, Saudi Arabia}
\icmlaffiliation{mipt}{Moscow Institute of Physics and Technology, Russia}
\icmlaffiliation{ed}{University of Edinburgh, United Kingdom}

\icmlcorrespondingauthor{Samuel Horv\'{a}th}{samuel.horvath@kaust.edu.sa}
\icmlcorrespondingauthor{Peter Richt\'{a}rik}{peter.richtarik@kaust.edu.sa}

\icmlkeywords{Optimization, Arbitrary Sampling, Importance Sampling}

\vskip 0.3in
]



\printAffiliationsAndNotice{\textsuperscript{*} This work was awarded \href{http://2018.ds3-datascience-polytechnique.fr/posters/}{\emph{the Best DS\textsuperscript{3} Poster Award} }  among 170 posters at the 2018  Data Science Summer School (DS\textsuperscript{3}) at Ecole Polytechnique in Paris, France. \\} 

\begin{abstract}
 We provide the first importance sampling variants of variance-reduced algorithms for empirical risk minimization with non-convex loss functions. In particular, we analyze non-convex versions of \texttt{SVRG}, \texttt{SAGA} and \texttt{SARAH}. Our methods have the capacity to speed up the training process by  an order of magnitude compared to the state of the art on real datasets. Moreover, we also improve upon current mini-batch analysis of these methods by proposing  importance sampling for minibatches in this setting. Ours are the first optimal samplings for minibatches in the literature on stochastic optimization. Surprisingly, our approach can in some regimes lead to superlinear speedup with respect to the minibatch size, which is not usually present in stochastic optimization. All the above results follow from a general analysis of the methods which works with {\em arbitrary sampling}, i.e., fully general randomized strategy for the selection of subsets of examples to be sampled in each iteration. Finally, we also perform a novel importance sampling analysis of \texttt{SARAH} in the convex setting.
\end{abstract}

\section{Introduction}

Empirical risk minimization (ERM) is a key problem in machine learning as it plays a key role in training supervised learning models, including classification and regression problems, such as support vector machine, logistic regression and deep learning.  A generic ERM problem has the finite-sum form
 \begin{equation}
\label{eq:1}
  \min_{x}\ f(x) \eqdef \frac{1}{n}\sum_{i=1}^n f_i(x),
\end{equation} 
where $x$ corresponds to the parameters defining a model, $f_i(x)$ is the loss of the model $x$ associated with  data point $i$, and $f$ is the average (empirical) loss across the entire training dataset. In this paper we  focus on the case when the functions $f_i$ are $L_i$--smooth but {\em non-convex}. We assume the problem has a solution $x^*$.

One of the most popular algorithms for solving \eqref{eq:1} is  stochastic gradient descent (\texttt{SGD}) \citep{nemirovsky1983problem, Nemirovski_2009}. In recent years, tremendous effort was exerted to improve its performance, leading to various enhancements which use acceleration  \cite{Katyusha}, momentum  \cite{BasicMomentum}, minibatching  \cite{pegasos2}, distributed implementation \cite{CoCoA+, CoCoA+OMS}, importance sampling  \cite{zhao2015stochastic, is-minibatches, Quartz, SCP}, higher-order information  \cite{SDNA, SBFGS}, and a number of other techniques. 

\subsection{Variance-reduced methods}

A particularly important recent advance has to do with the design of {\em variance-reduced (VR)} stochastic gradient methods, such as \texttt{SAG} \citep{SAG}, \texttt{SDCA} \citep{SDCA, UCDC}, \texttt{SVRG} \citep{SVRG}, \texttt{S2GD} \citep{S2GD},  \texttt{SAGA} \citep{SAGA}, \texttt{MISO} \citep{MISO}, \texttt{FINITO} \citep{FINITO} and \texttt{SARAH} \citep{SARAH},
 which  operate by modifying the classical stochastic gradient  direction  in each step of the training process in various clever ways so as to progressively reduce its variance as an estimator of the true gradient.  We note that \texttt{SAG} and \texttt{SARAH},  historically the oldest and one the newest VR methods in the list, respectively, use a biased estimator of the gradient. In theory, all these methods enjoy linear convergence rates on smooth and strongly convex functions, which is in contrast with slow sublinear rate of  \texttt{SGD}.  VR methods are also easier to implement as they do not rely on a decreasing learning rate schedule. VR methods were  recently extended to work with (non-strongly) convex losses by \citet{S2GD}, and more recently also to  non-convex losses by \citet{reddi2016fast, reddi2016stochastic, allen2016variance, SARAH-nonconvex} - in all cases leading to best current rates for  \eqref{eq:1} in a given function class. 

\subsection{Importance sampling, minibatching and non-convex models}

In the context of problem \eqref{eq:1}, importance sampling refers to the technique of assigning carefully designed {\em non-uniform} probabilities $\{p_i\}$ to the $n$ functions $\{f_i\}$, and using these, as opposed to uniform probabilities, to sample the next data point (stochastic gradient) during the training process.

 Despite the huge theoretical and practical success of VR methods, there are still considerable gaps in our understanding. For instance,  an importance sampling variant of the popular \texttt{SAGA} method, with the ``correct'' convergence rate, was only designed very recently by \citet{JacSketch}; and the analysis applies to strongly convex $f$ only. A coordinate descent variant of \texttt{SVRG} with importance sampling, also in the strongly convex case, was analyzed by \citet{S2CD}. However, the method does not seem to admit a fast implementation. For dual methods based on coordinate descent, importance sampling is relatively well understood \citep{RCDM, UCDC, ALPHA, Quartz, NUACDM}.

The territory is completely unmapped in the non-convex case, however. To the best of our knowledge, {\em no importance sampling} VR methods have been designed nor analyzed in the popular case when the functions $\{f_i\}$ are {\em non-convex}. An exception to this is  \texttt{dfSDCA} \citep{dfSDCA}; however, this method applies to an explicitly regularized version of \eqref{eq:1}, and while the individual functions are allowed to be non-convex, the average $f$ is assumed to be convex.  Given the dominance of stochastic gradient  type methods in training large non-convex models such as deep neural networks, theoretical investigation of VR methods that can benefit from importance sampling is much needed. 

The situation is worse still when one asks for {\em importance sampling of minibatches.} To the best of our knowledge, there are only a handful of papers on this topic \citep{NSync, is-minibatches, ACD2019}, none of which apply to the non-convex setting considered here, nor to the methods we will analyze, and the problem is open. This is despite the fact that minibatch methods are de-facto the norm for training  deep nets. In practice, typically relatively small ($\cO(1)$ or $\cO(\log n)$) minibatch sizes are used.
 
 \subsection{Contributions}

\begin{table*}[t]
\begin{center}
\begin{tabular}{|c||c|c|c|}
\hline
\bf Algorithm & \bf Uniform  sampling & \bf Arbitrary sampling {\bf [NEW]} &  \bf  $\Sam^*$ sampling {\bf [NEW]}\\
\hline 
\hline
& & & \\
\texttt{SVRG}    & $\max\left\{n, \frac{{\color{purple}(1+\nicefrac{4}{3})L_{\max}} c_1  n^{2/3} }{\epsilon} \right\}$  & $\max\left\{n,  \frac{ {\color{purple} (1+ \nicefrac{4\alpha}{3})\bar{L}} c_1 n^{2/3}}{\epsilon}\right\} $ & $\max \left\{n, \frac{ {\color{purple} \left(1+ \frac{4(n-b)}{3n}\right)\bar{L}} c_1 n^{2/3}}{\epsilon}\right\} $ \\
& & & \\
\texttt{SAGA}    & $n+ \frac{ {\color{purple} 2L_{\max}} c_2n^{2/3}}{\epsilon} $   & $ n+ \frac{{\color{purple}(1+\alpha) \bar{L}} c_2 n^{2/3}}{\epsilon} $ & $ n+ \frac{{\color{purple}(1+\frac{n-b}{n} )\bar{L}} c_2 n^{2/3}}{\epsilon} $ \\
& & & \\
\texttt{SARAH}  & $ n+ \frac{ {\color{purple}\frac{n-b}{n-1} L_{\max}^2} c_3}{\epsilon^2} $   & $ n+ \frac{ {\color{purple}\alpha\bar{L}^2}c_3}{\epsilon^2} $   & $ n+ \frac{ {\color{purple}\frac{n-b}{n}\bar{L}^2} c_3}{\epsilon^2}  $   \\ 
& & & \\
\hline
\end{tabular}

\end{center} 
\caption{Stochastic gradient evaluation complexity for achieving $\Exp{\| \nabla f(x) \|^2} \leq \epsilon$ for two variants of \texttt{SVRG}, \texttt{SAGA} and \texttt{SARAH} for minimizing the average of  smooth non-convex functions. Constants: $c_1, c_2, c_3$ are universal constant, $L_{\max} = \max_i L_i$; $\bar{L} = \tfrac{1}{n}\sum_i L_i$; $b$ = (average) minibatch size (hidden in $\alpha$); $\alpha$ can be for specific samplings smaller than $1$ and decreasing with  increasing $b$, which  can lead to superlinear speedup in $b$. For \texttt{SARAH} this guarantee holds for one outer loop with minibatch size, where we assume $16\bar{L}^2(f(x^0)-f(x^*)^2)/(\epsilon b)^2\gg 0$, in other words, minibatch size is not too big comparing to the required precision.}
\label{tab:algo-comparison}
\end{table*}

The main contributions of this paper are:

 {\bf Arbitrary sampling.} We peform a  general analysis of three popular VR methods---\texttt{SVRG} \citep{SVRG}, \texttt{SAGA} \citep{SAGA} and \texttt{SARAH} \citep{SARAH}---in the {\em arbitrary sampling} paradigm \citep{NSync,ALPHA, ESO, Quartz, SCP}. That is, we prove general complexity results which hold for an {\em arbitrary random set valued mapping} (aka arbitrary sampling) generating the minibatches of examples used by the algorithms in each iteration.  
 
{ \bf Optimal sampling.} Starting from our general complexity results which hold for arbitrary sampling (see the second column in Table~\ref{tab:algo-comparison}), we are able calculate the {\em optimal sampling out of all samplings of a given minibatch size} (see Lemma~\ref{optimality} and also the last column in Table~\ref{tab:algo-comparison}). {\em This is the first time an optimal minibatch sampling was computed (from the class of all samplings) in the literature for any stochastic optimization method we know, including all variants of \texttt{SGD} and coordinate descent.} Indeed, while the results in  \citep{NSync, is-minibatches, ACD2019} and other works on this topic construct importance sampling for minibatches, these are not shown nor believed to be optimal.

{\bf Improved rates.} Our iteration complexity bounds improve upon the best current rates for these methods even in the non-minibatch case. For \texttt{SVRG} and \texttt{SAGA}, this is true even when $L_i=L_j$ for all $i,j$, which is counter-intuitive as classical importance sampling is proportional to the constants $L_i$, which in this case would lead to uniform probabilities. Our importance sampling can be faster by up to the factor of $n$ compared to the current state of the art (see Table~\ref{tab:algo-comparison} and Appendix~\ref{sec:improvements}).  Our methods can enjoy {\em linear speedup} or even for some specific samplings  {\em superlinear speedup}  in minibatch size. That is, the number of iterations needed to output a solution of a given accuracy drops by a factor equal or greater to the minibatch size used. This is of utmost relevance to the practice of training  neural nets with minibatch stochastic methods as our results predict that this is to be expected. We design importance sampling and {\em approximate importance sampling for minibatches} which in our experiments vastly outperform the standard uniform minibatch strategies.

{\bf Best rates for SARAH under convexity.} Lastly, we also perform an analysis of importance sampling variant of \texttt{SARAH} in the convex and strongly convex case (Appendix~\ref{sec:SARAH-convex}). These are the currently fastest rates for \texttt{SARAH}.

\section{Importance Sampling for Minibatches}

As mentioned in the introduction, we assume throughout that $f_i:\R^d \mapsto \R$ are smooth, but not necessarily convex. In particular, we assume that $f_i$ is $L_i$--smooth; that is, 
$\| \nabla f_i(x) - \nabla f_i(y) \| \leq L_i \| x  - y \|$ for all $ x,y\in \R^d,$ where $\|x\|\eqdef (\sum_i x_i^2)^{1/2}$ is the standard Euclidean norm.
Let  us define $\bar{L} \eqdef \frac{1}{n}\sum_{i=1}^n L_i$. Without loss of generality assume that $L_1 \leq L_2 \leq \dots \leq L_n$.

In this work, our aim is to find an $\epsilon$--accurate solution in expectation.  A stochastic iterative algorithm for solving \eqref{eq:1} is said to achieve 
$\epsilon$--accurate solution if the random output $x_a$ of this algorithm satisfies $\Exp{\norm{\nabla f(x_a)}^2} \leq \epsilon.$

\subsection{Samplings} \label{sec:samplings}

Let $\Sam$ be a random set-valued mapping (``sampling'') with values in $2^{[n]}$, where $[n] \eqdef \{1,2,\dots,n\}$. A sampling is uniquely defined by assigning probabilities to all $2^n$ subsets of $[n]$. With each sampling we  associate a {\em probability matrix} $\mP \in \R^{n\times n}$ defined by \[\mP_{ij} \eqdef \Prob(\{i,j\}\subseteq \Sam).\] The {\em probability vector} associated with $\Sam$ is the vector composed of the diagonal entries of $\mP$: $p = (p_1,\dots,p_n)\in \R^n$, where $p_i\eqdef \Prob(i\in \Sam)$. We say that $\Sam$ is {\em proper} if $p_i>0$ for all $i$. It is easy to show that \begin{equation} \label{eq:b} b\eqdef \Exp{|\Sam|} = \trace{\mP} = \sum_{i=1}^n p_i.\end{equation}
From now on, we will refer to $b$ as the {\em minibatch size} of sampling $\Sam$.   
It is known that $\mP - p p^\top\succeq 0$  \citep{PCDM}. Let us without loss of generality assume that $p_1 \leq p_2 \leq \dots \leq p_n $ and define constant $k = k(S) \eqdef |\{ i \in [n]: p_i < 1\}| = \max\{i: p_i < 1\}$ to be the number of $p_i$'s not equal to one.

While our complexity results are general in the sense that they hold for any proper sampling, we shall now consider three special samplings; all with minibatch size $b \in (0,n]$:

\paragraph{1. Standard uniform minibatch sampling ($\Sam=\Sam^u$).}  $\Sam$ is chosen uniformly at random from all subsets of $[n]$ of cardinality $b$. Clearly, $|\Sam|=b$ with probability 1. The probability matrix is given by \begin{equation}\label{eq:prob_matrix_nice}\mP_{ij} = \begin{cases} \tfrac{b}{n} & i=j,\\
\tfrac{b(b-1)}{n(n-1)} & i\neq j .
\end{cases}\end{equation}
In the literature, this is known as the $b$--nice sampling \citep{PCDM, ESO}.

\paragraph{2. Independent sampling $(\Sam=\Sam^*)$.} For each $i\in [n]$ we independently flip a coin, and with probability $p_i>0$ include element $i$ into $\Sam$. Hence, by construction, $p_i = \Prob(i\in \Sam)$ and $\Exp{|\Sam|} \overset{\eqref{eq:b}}{=} \sum_{i} p_i = b$. The probability matrix of $\Sam$ is
 \[\mP_{ij} = \begin{cases} p_i & i=j,\\
p_i p_j & i\neq j .
\end{cases}\]

\paragraph{3. Approximate independent sampling $(\Sam=\Sam^a)$.} Independent sampling  has the disadvantage that $k=k(S)$ coin tosses need to be performed in order to generate the random set.  However, we would like to sample at the cost $\cO(b+k-n)$ coin tosses instead. We now design a sampling which has this property and which in a certain precise sense, as we shall see later, approximates the independent sampling.
In particular, given an independent sampling with parameters $p_i$ for $i\in [n]$, let $a = \lceil k \max_{i\leq k} p_i \rceil$. Since $\max_{i \leq k} p_i \geq \tfrac{b+k-n}{k}$, it follows that $a \geq b+k-n$. On the other hand, if $\max_{i \leq k} p_i=\cO((b+k-n)/k)$, then $a=\cO(b+k-n)$. We now sample a single set $\Sam'$ of cardinality $a$ using the standard uniform minibatch sampling (just for $i \leq k$). Subsequently, we apply an independent sampling to select elements of $\Sam'$, with selection probabilities $p'_i = k p_i/a$. The resulting set is $\Sam$. Since \[\Prob(i \in \Sam) =  \tfrac{\binom{k}{a-1}}{\binom{k}{a}} \frac{kp_i}{a}  = p_i,\]  \[\Prob(\{i,j\} \subseteq \Sam) =  \tfrac{\binom{k-2}{a-2}}{\binom{k}{a}}\frac{kp_i}{a} \frac{np_j}{a}   = \frac{(a-1)k}{a(k-1)}p_i p_j \] for $i,j \leq k$, the probability matrix of $\Sam$ is given by 
 \[\mP_{ij} = \begin{cases} p_i & i=j, \\
 \tfrac{(a-1)k}{a(k-1)} p_i p_j & i\neq j; i,j\leq k, \\
 p_i p_j & \text{otherwise} .
\end{cases}\]
Since $ \tfrac{(a-1)k}{a(k-1)} \approx 1$, the probability matrix of the approximate independent sampling approximates that of the independent sampling. Note that $\Sam$ includes both the standard uniform minibatch sampling and the independent sampling as special cases. Indeed, the former is obtained by choosing $p_i=b/n$ for all $i$ (whence $a=b$ and $p'_i = 1$ for all $i$), and the latter is obtained by choosing $a=n$ instead of $a=\lceil k \max_{i\leq k} p_i \rceil$.

\subsection{Key lemma}

The following lemma, which we use as an upper bound for variance,  plays a key role in our analysis.

\begin{lemma}
\label{lem:upperv}
Let $\zeta_1, \zeta_2,\dots,\zeta_n$ be vectors in $\R^d$ and let $\bar{\zeta}\eqdef\frac1n \sum^n_{i=1}\zeta_i$ be their average. Let $\Sam$ be a proper sampling (i.e., assume that $p_i=\Prob(i\in \Sam)>0$ for all $i$).  Assume that there is $v \in \R^n$ such that 
\begin{equation}
\label{eq:ESO}
\mP - pp^\top \preceq \Diag{p_1 v_1, p_2 v_2, \dots, p_n v_n}.
\end{equation}
Then 
\begin{equation}\label{eq:key_inequality}
 \E{\left\|   \sum_{i \in \Sam} \frac{\zeta_i}{np_i} - \bar{\zeta}\right\|^2} \leq  \frac{1}{n^2} \sum_{i=1}^n \frac{v_i}{p_i}  \|\zeta_i\|^2,
\end{equation}
where the expectation is taken over sampling $\Sam$. Whenever \eqref{eq:ESO} holds, it must be the case that 
\begin{equation}
\label{v_i}
 v_i \geq 1 - p_i.
\end{equation}
Moreover, \eqref{eq:ESO} is always satisfied for  $v_i=n(1-p_i)$ for $i \leq k$ and $0$ otherwise. Further, if  $|\Sam| \leq d$ with probability 1 for some $d$, then \eqref{eq:ESO} holds for $v_i = d$. The standard uniform minibatch sampling admits $v_i = \frac{n-b}{n-1}$, the independent sampling  admits $v_i = 1 - p_i$, and the approximate independent sampling admits the choice \[v_i = 1 - p_i(1 - \tfrac{k-a}{a(k-1)})\] if $i \leq k$, $v_i=0$ otherwise.
\end{lemma}

\subsection{Optimal sampling}

The following quantities play a key role in our general complexity results:
\begin{equation}
\label{def_K} 
K \eqdef \frac{b}{n^2} \sum_{i=1}^n \frac{v_i L_i^2}{p_i}, \qquad \alpha \eqdef \frac{K}{\bar{L}^2}.
\end{equation}  
Above, $b=\Exp{|S|}$ is the minibatch size, $p_i=\Prob(i\in S)$, $\bar{L}\eqdef \tfrac{1}{n}\sum_i L_i$ and $\{v_i\}$ are defined in \eqref{eq:ESO} in Lemma~\ref{lem:upperv}.

Our theory shows (see the 2nd column of Table~\ref{tab:algo-comparison} for a summary, and Section~\ref{sec:rates} for the full results) that in order to optimize the iteration complexity, we need to design sampling $\Sam$ for which the value $\alpha$ is as small as  possible. The following result sheds light on how $\Sam$ should be chosen, from samplings of a given minibatch size $b$, to minimize $\alpha$. 

\begin{lemma}
\label{optimality}
Fix a minibatch size $b \in (0,n]$.  Then the quantity $\alpha$, defined in \eqref{def_K}, is minimized for the choice $\Sam=\Sam^*$ with the probabilities
\begin{equation}
\label{prob}
p_i \eqdef  \begin{cases}
       (b +  k - n)\frac{L_i}{\sum_{j=1}^k L_j,} &\quad\text{if } i \leq k  \\
       1, &\quad\text{if } i > k
\end{cases}, 
\end{equation}
where $k$ is the largest integer satisfying $0 < b + k - n \leq \nicefrac{\sum_{i =1}^k L_i}{L_k}$ (for instance, $k=n-b+1$ satisfies this).  Usually, if $L_i$'s are not too much different, then $k = \cO(n)$, for instance, if $b L_n \leq \sum_{i=1}^nL_i $ then $k = n$.
If we choose $\Sam=\Sam^a$, then $\alpha$ is minimized for \eqref{prob} with 
\begin{equation}
\label{alpha_im}
\alpha = \left(\frac{b \left(\sum_{i=1}^k L_i\right)^2}{(b + k - n)n^2} - \frac{bs}{n^2} \sum_{i=1}^k L_i^2\right) / \bar{L}^2,
\end{equation}
where $s = 1$ for $\Sam^*$  and $s = 1 - \tfrac{k-a}{a(k-1)}$ for $\Sam^a$.  Moreover, if we assume\footnote{Note, that this can be always satisfied, if we uplift the smallest $L_i$'s, because if function is $L$-smooth, then it is also smooth with $L'\geq L$.}  $b L_n \leq \sum_{i = 1}^n L_i$, then $k = n$, thus
\begin{eqnarray}
\alpha_{S^*} &=& 1 - b \frac{\sum_{i=1}^n L_i^2}{\left(\sum_{i=1}^n L_i\right)^2} \leq \frac{n-b}{n}, \notag \\
\alpha_{S^u}&=& (n-b)\frac{n}{n-1} \frac{\sum_{i=1}^n L_i^2}{\left( \sum_{i=1}^n L_i\right)^2}. \notag
\end{eqnarray}
For the special case of all $L_i's$ are the same, one obtains
\[
\alpha_{S^*} =  \frac{n-b}{n}, \quad  \alpha_{S^u}= \frac{n-b}{n-1}. \]
\end{lemma} 

From now on, let $S^*, S^a$ denote \textit{Independent Sampling} and  \textit{Approximate Inpedendent Sampling}, respectively, with the probabilities defined in \eqref{prob}. Lemma~\ref{optimality} guarantees that the sampling  $\Sam^*$ is optimal (i.e., minimizes $\alpha$). Moreover, if we let $b_{\max} \eqdef  \max \left\{b \;|\; bL_n \leq \sum_i L_i \right\},$ then we obtain {\em superlinear speedup in $b$}, up to $b_{\max}$ for all three algorithms.

\section{SVRG, SAGA and SARAH}
\label{sec:rates}

In all of the results of this section we assume that $\Sam$ is an {\em arbitrary proper sampling}. Let $b=\E{|\Sam|}$ be the (average) minibatch size. We assume that $v$ satisfies \eqref{eq:ESO} and that $\alpha$ (which depends on $v$) is defined as in \eqref{def_K}. All complexity results will depend  on $\alpha$ and $b$. 

We propose three methods, Algorithm~\ref{alg:minibatch-svrg}, \ref{alg:minibatch-saga} and \ref{alg:minibatch-SARAH}, which are generalizations  of 
original \texttt{SVRG} \cite{reddi2016stochastic}, \texttt{SAGA}  \cite{reddi2016fast} and \texttt{SARAH} \cite{SARAH-nonconvex} to the arbitrary sampling setting, respectively. The original non-minibatch methods arise as special cases for the sampling  $\Sam = \{i\}$ with probability $1/n$, and the original minibatch methods arise as a special case for the sampling $\Sam^u$ (described in Section~\ref{sec:samplings}).

\label{mb_SVRG}
\begin{algorithm}[tb]
   \caption{\texttt{SVRG} with arb.\ sampling $\left(x^0,m,T,\eta, \Sam\right)$}
   \label{alg:minibatch-svrg}
\begin{algorithmic}[1]
   \STATE $\tilde{x}^0 = x^0_m = x^0 $, $M = \lceil T/m \rceil$
   \FOR{$s=0$ {\bfseries to} $M-1$}
   \STATE $x^{s+1}_0 = x^{s}_m$
   \STATE $g^{s+1} = \frac{1}{n} \sum_{i=1}^n \nabla f_{i}(\tilde{x}^{s})$
   \FOR{$t=0$ {\bfseries to} $m-1$}
   \STATE Draw a random subset  (minibatch) $\Sam_t  \sim \Sam$
   \STATE $\displaystyle v_t^{s+1} =   \sum_{i_t \in \Sam_t} \tfrac{\nabla f_{i_t}(x^{s+1}_t) - \nabla f_{i_t}(\tilde{x}^{s})}{np_{i_t}}+ g^{s+1}$
   \STATE $x^{s+1}_{t+1} = x^{s+1}_{t} - \eta v_t^{s+1} $
   \ENDFOR
   \STATE $\tilde{x}^{s+1} = x_{m}^{s+1}$
   \ENDFOR
   \STATE {\bfseries Output:} Iterate $x_a$ chosen uniformly at random from $\{\{x^{s+1}_t\}_{t=0}^{m}\}_{s=0}^{M}$.
\end{algorithmic}
\end{algorithm}

Our general result for \texttt{SVRG} follows. 
 
\begin{theorem}[Complexity of \texttt{SVRG} with arbitrary sampling]
There exist universal constants $\mu_2>0$, $0<\nu_2 < 1$ such that the output of Alg.~\ref{alg:minibatch-svrg} with mini-batch size $b \leq \alpha n^{2/3}$, step size $\eta = \mu_2 b/(\alpha \bar{L}n^{2/3})$, and parameters $\beta = \bar{L}/n^{1/3}$, $m = \lfloor n\alpha/(3b\mu_2) \rfloor$ and $T$ (multiple of $m$) satisfies:
\[
\Exp{\|\nabla f(x_a)\|^2} \leq \frac{\alpha\bar{L}n^{2/3} [f(x^{0}) - f(x^*)]}{bT\nu_2} . 
\] 
Thus in terms of stochastic gradient evaluations to obtain $\epsilon$-accurate solution, one needs following number of iterations
\begin{equation}
\max\left\{n,\frac{\mu_2 \bar{L}n^{(2/3)}(f(x^0)-f(x^*))}{\epsilon\nu_2}\left(1+\frac{\alpha}{3\mu_2}\right)\right\}. \notag
\end{equation}  
\label{thm:nonconvex-minibatch}
\end{theorem}

In the next theorem we provide a generalization of  the results by \citet{reddi2016fast}.

\begin{algorithm}[tb]
   \caption{\texttt{SAGA} with arbitrary sampling $\left(x^0,d,T,\eta, \Sam\right)$}
   \label{alg:minibatch-saga}
\begin{algorithmic}[1]
   \STATE  $\alpha_{i}^0 = x^0$ for $i \in [n]$, 
   \STATE $g^{0} = \frac{1}{n} \sum_{i=1}^n \nabla f_{i}(\alpha_i^{0})$
   \FOR{$t=0$ {\bfseries to} $T-1$}
   \STATE  Draw a random subset (minibatch) $\Sam_t \sim \Sam$
   \STATE  Pick  random subset $J_t\subset [n]$ s.t. $\Prob(j \in J_t) = \tfrac{d}{n}$ 
   \STATE $\displaystyle v^t =  \sum_{i \in \Sam_t} \tfrac{\nabla f_{i}(x^t) - \nabla f_{i}(\alpha_{i_t}^{t})}{np_{i_t}} + g^{t}$
   \STATE $x^{t+1} = x^{t} - \eta v^t$
   \STATE $\alpha_{j}^{t+1} = x^t$ for $j \in J_t$ and $\alpha_{j}^{t+1} = \alpha_{j}^{t}$ for $j \notin J_t$
   \STATE $g^{t+1} = g^t - \frac{1}{n} \sum_{j \in J_t}(\nabla f_{j}(\alpha_{j}^t) - \nabla f_{j}(\alpha_{j}^{t+1}))$
   \ENDFOR
   \STATE {\bfseries Output:} Iterate $x_a$ chosen uniformly at random from $\{x^t\}_{t=0}^{T}$.
\end{algorithmic}
\end{algorithm}

\begin{theorem}[Complexity of \texttt{SAGA} with arbitrary sampling]
\label{thm:nonconvex-minibatch_SAGA}
There exist universal constants $\mu_3>0$, $0<\nu_3 < 1$ such that the output of Alg.~\ref{alg:minibatch-saga} with mini-batch size $b \leq \alpha n^{2/3}$, step size $\eta = b/(\mu_3\alpha \bar{L}^2 n^{2/3})$, and parameter $d = b/\alpha$  satisfies:
\[
\Exp{\|\nabla f(x_a)\|^2} \leq \frac{\alpha\bar{L}n^{2/3} [f(x^{0}) - f(x^*)]}{bT\nu_3} .
\]
Thus, in terms of stochastic gradient evaluations, to obtain $\epsilon$--accurate solution, one needs following number of iterations
\begin{equation}
 n + \frac{ \bar{L}n^{(2/3)}(f(x^0)-f(x^*))}{\epsilon\nu_3}(1+\alpha). \notag	
\end{equation} 
\end{theorem}


We  now introduce Algorithm~\ref{alg:minibatch-SARAH}: a general form of the \texttt{SARAH} algorithm  \citep{SARAH-nonconvex}. 


 \begin{algorithm}[tb]
   \caption{\texttt{SARAH} with arb.\ sampling $\left(x^0,m,T, \eta, \Sam\right)$}
   \label{alg:minibatch-SARAH}
 \begin{algorithmic}[1]
   \STATE $x_{0}^{0} = x^0$
    \FOR{$s=1$ {\bfseries to} $M-1$}
   \STATE $v^0_s = \frac{1}{n}\sum_{i=1}^{n} \nabla f_i(x^0_s)$
   \STATE  $x^1 = x^0 - \eta v^0$
   \FOR{$t=1$ {\bfseries to} $m-1$}
   \STATE  Draw a random subset (minibatch) $\Sam_t \sim \Sam$
   \STATE $v^{t} = \sum_{i \in \Sam_t}\frac{1}{np_i} (\nabla f_{i} (x^{t}) - \nabla f_{i}(x^{t-1})) + v^{t-1}$
   \STATE $x^{t+1} = x^{t} - \eta v^{t}$
   \ENDFOR
   \STATE $x^0_{s+1}$ chosen uniformly at randomly from $\{x^t_s\}_{t=0}^{m}$
   \ENDFOR 
   \STATE {\bfseries Output:} Iterate $x_a = x^0_{M}$ 
\end{algorithmic}
\end{algorithm}

\begin{theorem}[Complexity of \texttt{SARAH} with arbitrary sampling]\label{thm:nonconvex_01_mb}
Consider one outer loop of Alg.~\ref{alg:minibatch-SARAH} with 
\begin{eqnarray}
\label{eta_mb}
\eta \leq \tfrac{2}{\bar{L}\left(\sqrt{1 + \frac{4\alpha m}{b}} + 1\right)}.
\end{eqnarray}
Then the output $x_a$ satisfies:
\[
\Exp{ \| \nabla f(x_a)\|^2 } \leq \frac{2}{\eta(m+1)} [ f(x^{0}_s) - f(x^{*})].\]
Thus, to obtain $\epsilon$--accurate solution, one needs
\begin{eqnarray}
&n + \frac{16 \alpha\bar{L}^2 (f(x^0)-f(x^*))^2}{2\epsilon^2} + \notag \\
&\frac{\sqrt{16^2 \alpha^2\bar{L}^4 (f(x^0)-f(x^*))^4 + 16 \epsilon^2\bar{L}^2(f(x^0)-f(x^*))^2b^2}}{2\epsilon^2}
\notag
\end{eqnarray}
stochastic gradient evaluations. 
\end{theorem}

If all $L_i$'s are the same and we choose $\Sam$ to be $\Sam^a$, thus uniform with mini-batch size $b$, we can get back original result from \cite{SARAH-nonconvex}. Taking $b = n$, we can restore gradient descent with the correct step size $\nicefrac{1}{\bar{L}}$. 

\section{Additional Results}

In this section we describe three additional results: linear convergence for {\tt SVRG}, {\tt SAGA} and {\tt SARAH} for gradient dominated functions, the first importance sampling results for {\tt SARAH} for convex functions, and an array of new rates (with slight improvements) for non-minibatch versions of the above three methods for non-convex problems.

\subsection{Gradient dominated functions}

\begin{definition}
\label{gd-dom}
We say that $f$ is $\tau$-gradient dominated if  
\[
f(x) - f(x^*) \leq \tau \| \nabla f(x) \|^2,
\]
for all $x \in \R^d$, where $x^*$ is an optimal solution of  \eqref{eq:1}.
\end{definition}

\begin{algorithm}[tb]
   \caption{GD-Algorithm$\left(x^0,T, \mathbb{A} \right)$}
   \label{alg:gd-alg}
\begin{algorithmic}
   \STATE {\bfseries Input:} $x^0 \in \R^d$, $T$, $\mathbb{A}$
   \FOR{$k=0$ to $K$}
   \STATE $x^{k} = \text{Non-convex algorithm}(x^{k-1}, T, \mathbb{A})$
   \ENDFOR
   \STATE {\bfseries Output:} $x^K$
\end{algorithmic}
\end{algorithm}

Gradient dominance is a weaker version of strong convexity due to the fact that if function is $\mu$-strongly convex then it is $\tau$-gradient dominated, where $\tau = 2/\mu$.

Any of the non-convex methods in this paper  can be used as a subroutine of  Algorithm~\ref{alg:gd-alg}, where $T$ is the number of steps of the subroutine and $\mathbb{A}$ is the set of optimal  parameters for the subroutine. We set $T= \alpha n^{2/3}/(b \nu_2)$ for \texttt{SVRG} and $T = \alpha n^{2/3}/(b\nu_3)$ for \texttt{SAGA}. In the case of \texttt{SARAH}, $T$ is obtained by solving $m+1 = 2/\eta$ in $m$ and setting $T\leftarrow m$. Using  Theorems~\ref{thm:nonconvex-minibatch}, \ref{thm:nonconvex-minibatch_SAGA}, \ref{thm:nonconvex_01_mb} and the above special choice of $T$,
we get
\begin{equation}
\Exp{ \| \nabla f(x^k) \|^2} \leq \frac{1}{2 \tau} \left(\Exp{f(x^{k-1})} - f(x^*)\right). \notag
\end{equation}
Combined with Definition~\ref{gd-dom}, this guarantees a linear convergence with the same constant terms consisting of $\alpha, \bar{L}$ and $b$ that we had before in our analysis.  

\subsection{Importance sampling for SARAH under convexity}

In addition to the results presented in previous sections, we also establish  importance sampling results for \texttt{SARAH} in convex and strongly convex cases (Appendix~\ref{sec:SARAH-convex}) with similar improvements as for the non-convex algorithm. Ours are the best current rates for \texttt{SARAH} in these settings.
 
 \subsection{Better rates for non-minibatch methods for non-convex problems}
 
Lastly, we also provide specialized non-minibatch versions of non-convex \texttt{SAGA}, \texttt{SARAH} and \texttt{SVRG}, which are special cases of their minibatch versions presented in the main part with slightly improved guarantees (see Theorems~\ref{thm:nonconvex-inter_is}, \ref{thm:nonconvex-gen_is},  \ref{thm:nonconvex-inter} and  \ref{thm:nonconvex-gen} in the Appendix).

\begin{figure*}[t]
\begin{multicols}{3}
\includegraphics[width = 0.33\textwidth ]{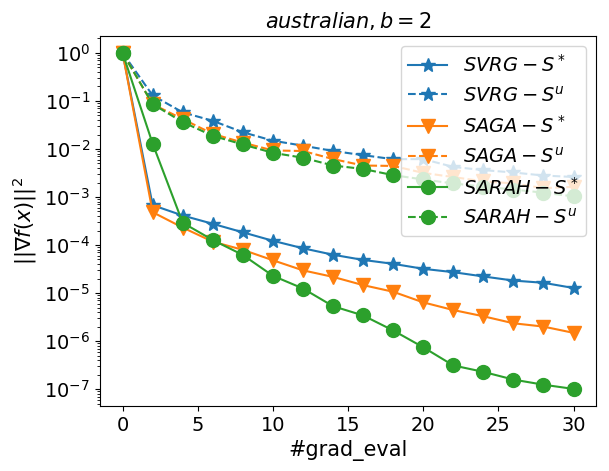}\par
\includegraphics[width =  0.33\textwidth ]{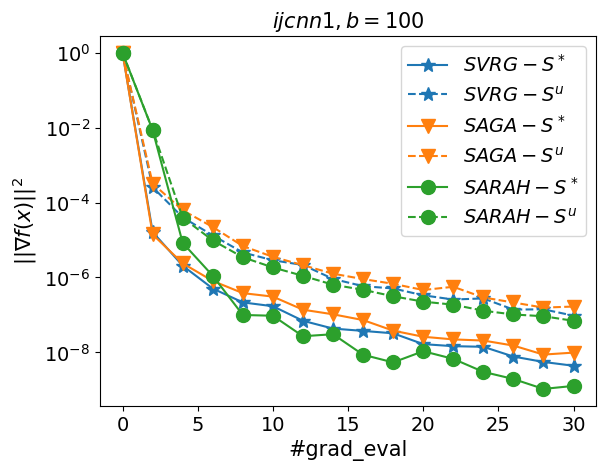}\par
\includegraphics[width =  0.33\textwidth ]{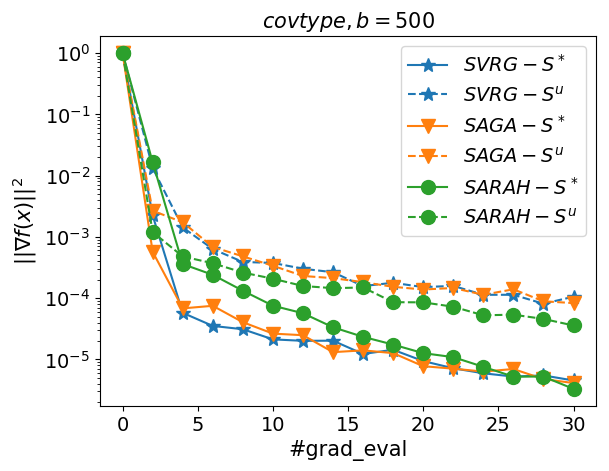}
\end{multicols}
\caption{Comparison of all methods with uniform and importance sampling: gradient norm.}
\label{fig:comp_all}
\end{figure*}

\section{Experiments}

In this section, we perform  experiments with regression for binary classification, where our loss function has the form 
\[f(x) = \frac{1}{n}\sum_{i=1}^n (1-y_i \sigma(a_i^\top x))^2,\] where $\sigma(z)$ is the sigmoid function. Hence, $f$  is smooth but non-convex. We chose this function because $L_i$'s can be easily computed (however,  in many cases, even for much more complex problems, they can usually be estimated). We use four LIBSVM datasets\footnote{The LIBSVM dataset collection is available at \url{https://www.csie.ntu.edu.tw/~cjlin/libsvmtools/datasets/}}: \textit{covtype, ijcnn1, splice, australian}. 


 Parameters of each algorithm are chosen as suggested by  the theorems in  Section~\ref{sec:rates}, and $x^0=0$. For \texttt{SARAH}, we chose $m=\lceil n/b \rceil $. The $y$ axis in all plots displays the norm of the gradient ($\|\nabla f(x)\|^2$) or the function value $f(x)$, and the $x$ axis depicts either epochs (1 epoch = 1 pass over data) or iterations.
 
%
%
%

\subsection{Importance vs uniform sampling}

Here we provide comparison of the methods with uniform ($S^u$) and  importance ($S^*$) sampling. Looking at  Figure~\ref{fig:comp_all}, one can see that  importance sampling outperforms uniform sampling for all three methods, in some cases, even by {\em several orders of magnitude.} For instance, in  the left plot (minibatch size $b=2$ and {\em australian} dataset) the improvement is as large as  $4$ orders of magnitude.

\begin{figure}[h]
\centering
\includegraphics[width = 0.38\textwidth ]{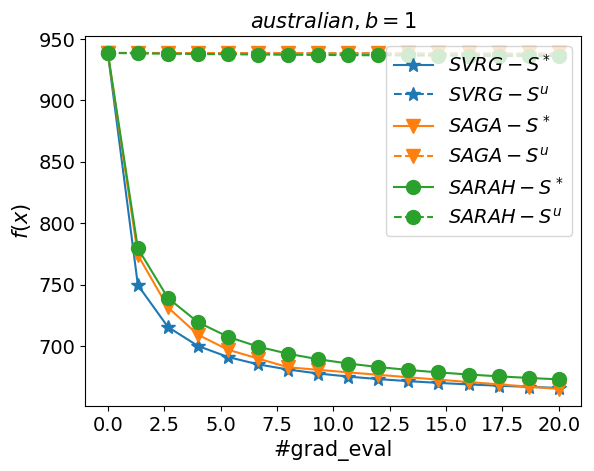}
\includegraphics[width = 0.38\textwidth ]{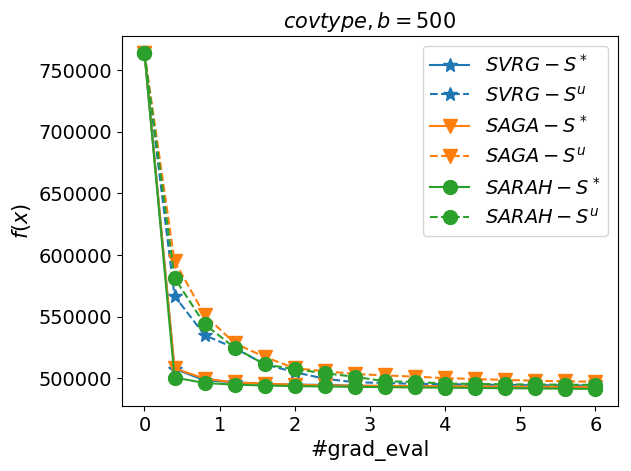}
\caption{Comparison of all methods with uniform and importance sampling for function values.}
\label{fig:comp_all_values}
\end{figure}

Looking at Figure~\ref{fig:comp_all_values}, one can see that there is an improvement not just in the norm of the gradient, but also in the function value. In the case of the \textit{australian} dataset, the  constants $L_i$'s are very non-uniform, and we can see that  the improvement is very significant.

\subsection{Linear or superlinear speedup}

Our theory suggests that linear or even superlinear speedup (in minibatch size $b$) can be obtained using the optimal independent $S^*$. Our experiments show that this is indeed the case in practice as well, and for all three algorithms.

\begin{figure*}[ht]
\begin{multicols}{3}
\includegraphics[width = 0.30\textwidth ]{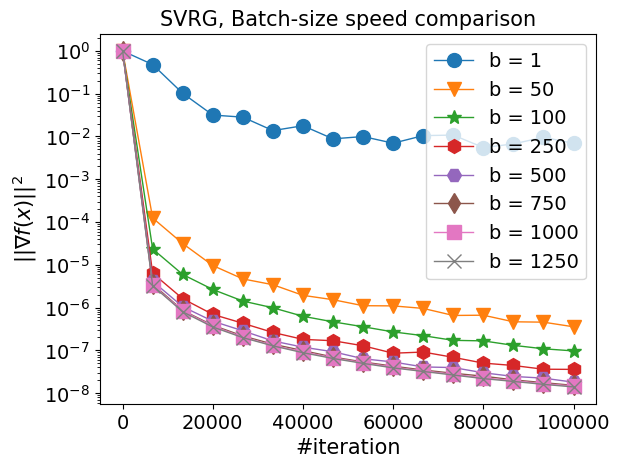}\par
\includegraphics[width = 0.30\textwidth ]{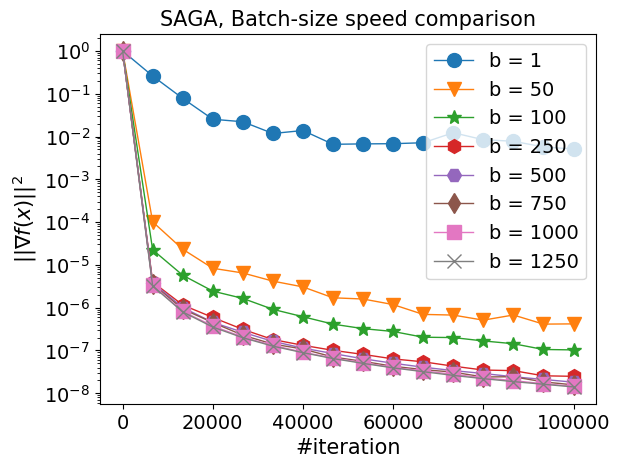}\par
\includegraphics[width =  0.30\textwidth ]{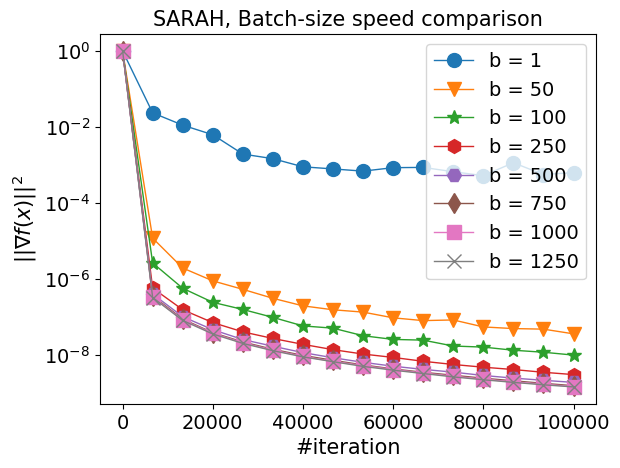}
\end{multicols}
\begin{multicols}{3}
\includegraphics[width =  0.30\textwidth ]{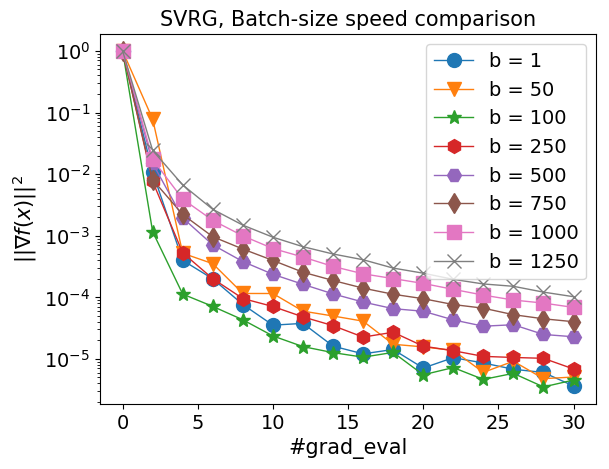}\par
\includegraphics[width =  0.30\textwidth ]{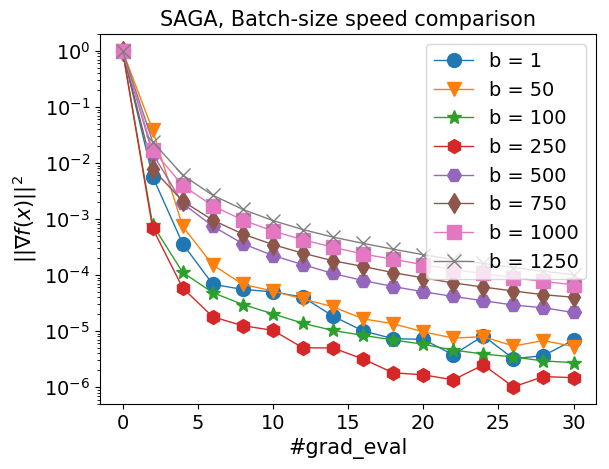}\par
\includegraphics[width =  0.30\textwidth ]{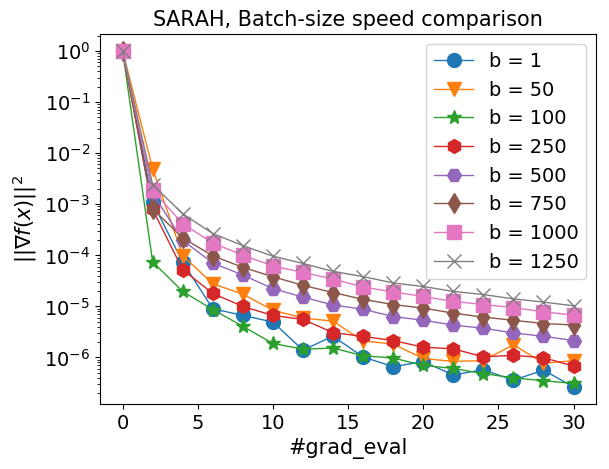}
\end{multicols}
\caption{Minibatch speedup, \textit{ijcnn1} dataset}
\label{fig:speed_up}
\begin{multicols}{3}
\includegraphics[width =  0.30\textwidth ]{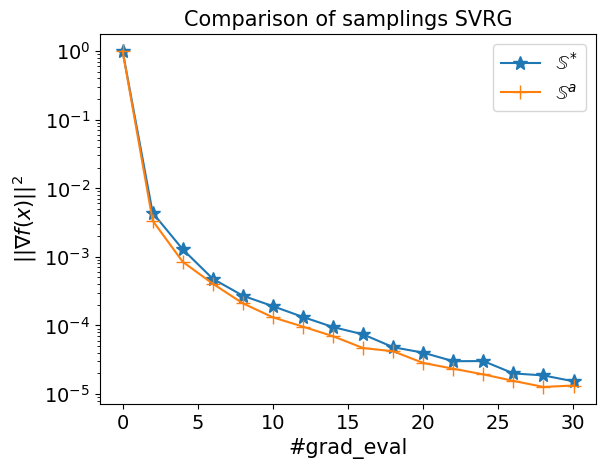}\par
\includegraphics[width =  0.30\textwidth ]{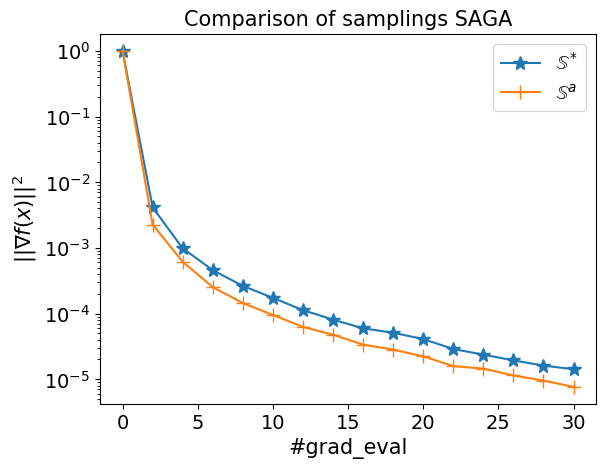}\par
\includegraphics[width =  0.30\textwidth ]{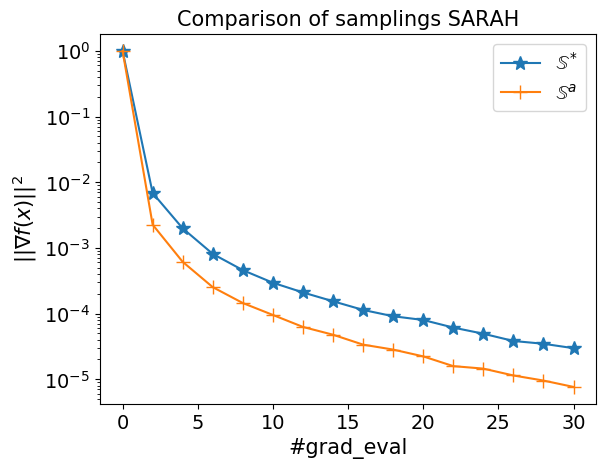}
\end{multicols}
\caption{Performance of sampling $\Sam^*$ vs. $\Sam^a$, \textit{splice} dataset.}
\label{fig:Samplings_comp}
\end{figure*}

Figure~\ref{fig:speed_up} confirms that linear, and sometimes even superlinear, speedup  is present. For this dataset, such speedup is present  up to the minibatch size of $250$. The plots in the top row of Figure~\ref{fig:speed_up} depict  convergence in a  simulated multi-core setting, where the number of cores is the same as the minibatch size.


\subsection{Independent  vs approximate independent sampling}

According to our theory, \textit{Independent Sampling} $\Sam^*$ is slightly better than \textit{Approximate Independent Sampling} $\Sam^a$. However, it is more expensive to use it in practice as generating samples from it involves more  computational effort for large $n$.  

The goal of our next experiment is to show that in practice $\Sam^a$ yields comparable or even faster convergence. Hence,  it is more reasonable to use this sampling for datasets where  the number of data points $n$ is big. For an efficient  implementation of $\Sam^a$,  we  can almost get rid of dependence on $n$. Intuitively, $\Sam^a$ works better because it has smaller variance in minibatch size than $\Sam^*$.

Indeed, it can be seen from Figure~\ref{fig:Samplings_comp} that $\Sam^a$ can outperform $\Sam^*$ in practice even though $\Sam^*$ is optimal in theory. The difference can be small (the left and the middle plot), but also quite significant (right plot).

\bibliographystyle{plainnat} 
\bibliography{literature}

\onecolumn
\clearpage
\appendix
\part*{Supplementary Material}
\setlength{\footskip}{20pt}

\section{Proof of Lemma~\ref{lem:upperv}}
\begin{proof}
Let $1_{i\in \Sam} = 1$ if $i\in \Sam$ and $1_{i\in \Sam} = 0$ otherwise. Likewise, let $1_{i,j\in \Sam} = 1$ if $i,j\in \Sam$ and  $1_{i,j\in \Sam} = 0$ otherwise. Note that $\Exp{1_{i\in \Sam}} = p_i$ and $\Exp{1_{i,j\in \Sam}}=p_{ij}$. Next, let us compute the mean of $X\eqdef  \sum_{i\in \Sam}\frac{\zeta_i}{np_i}$:
\begin{equation}\label{eq:098h0hf09h}\Exp{X} = \Exp{ \sum_{i\in \Sam}\frac{\zeta_i}{np_i}} = \Exp{\sum_{i=1}^n \frac{\zeta_i}{np_i} 1_{i\in \Sam}} = \sum_{i=1}^n \frac{\zeta_i}{np_i} \Exp{1_{i\in \Sam}}   = \frac{1}{n} \sum_{i=1}^n \zeta_i = \bar{\zeta}.\end{equation}

Let $\mA=[a_1,\dots,a_n]\in \R^{d\times n}$, where $a_i=\frac{\zeta_i}{p_i}$, and let $e$ be the vector of all ones in $\R^n$. We now write the  variance of $X$ in a  form which will be convenient to establish a bound:
\begin{eqnarray}
\Exp{\left\|X - \Exp{X}\right\|^2} &=&  \Exp{\norm{X}^2 } - \|\Exp{X}\|^2\notag \\
 &= & \Exp { \left\|   \sum_{i \in \Sam} \frac{\zeta_i}{np_i} \right\|^2 } - \norm{\bar{\zeta}}^2  \notag \\
&= &\Exp{ \sum_{i,j} \frac{\zeta_i^\top}{np_i } \frac{\zeta_j}{np_j  } 1_{i,j \in \Sam}}  - \norm{\bar{\zeta}}^2 \notag \\
&= &\sum_{i , j}  p_{ij}\frac{\zeta_i^\top}{np_i} \frac{\zeta_j}{np_j} - \sum_{i,j} \frac{\zeta_i^\top}{n} \frac{\zeta_j}{n} \notag \\
&=&\frac{1}{n^2} \sum_{i,j} (p_{ij} - p_i p_j) a_i^\top a_j  \notag \\
&=& \frac{1}{n^2} e^\top \left(\left( \mP - p p^\top \right) \circ \mA^\top \mA\right)e.\label{eq:n98dhg8f}
\end{eqnarray}
Since by assumption we have $\mP-pp^\top\preceq \Diag{p\circ v}$, we can further bound
\[
e^\top \left(\left( \mP - p p^\top \right) \circ \mA^\top \mA\right) e  \leq  e^\top \left( \Diag{p\circ v} \circ \mA^\top \mA \right) e  = \sum_{i=1}^n p_i v_i \|a_i\|^2.
\]
To obtain \eqref{eq:key_inequality}, it remains to combine this with \eqref{eq:n98dhg8f}. 

Inequality \eqref{v_i} follows by comparing the diagonal elements of the two matrices in \eqref{eq:ESO}.  Let us now verify the formulas for $v$. 
\begin{itemize}
\item Since $\mP-pp^\top$ is positive semidefinite \citep{PCDM}, we can bound $\mP-pp^\top \preceq n \Diag{\mP-p p^\top} = \Diag{p\circ v}$, where $v_i = n(1-p_i)$. 
\item It was shown by \citet[Theorem 4.1]{ESO} that $\mP \preceq d \Diag{p}$ provided that $|\Sam| \leq d$ with probability 1. Hence, $\mP - pp^\top \preceq \mP \preceq d \Diag{p}$, which means that $v_i=d$ for all $i$.
\item Consider now the independent sampling. Clearly,
\begin{eqnarray}
\mP - pp^\top = \begin{bmatrix}
    p_{1}(1-p_1) & 0& \dots  & 0 \\
    0 & p_{2}(1-p_2) & \dots  & 0 \\
    \vdots & \vdots & \ddots & \vdots \\
    0 & 0  & \dots  & p_{n}(1-p_n)
    \end{bmatrix} = 
    \Diag{p_1 v_1,\dots, p_n v_n}, \notag
\end{eqnarray}
where $v_i=  1-p_i$.

\item Consider the $b$--nice sampling (standard uniform minibatch sampling).  Direct computation shows that the probability matrix is given by 
\begin{eqnarray}\mP  = \begin{bmatrix}
    \tfrac{b}{n} & \tfrac{b(b-1)}{n(n-1)}& \dots  & \tfrac{b(b-1)}{n(n-1)} \\
    \tfrac{b(b-1)}{n(n-1)} & \tfrac{b}{n} & \dots  & \tfrac{b(b-1)}{n(n-1)} \\
    \vdots & \vdots & \ddots & \vdots \\
    \tfrac{b(b-1)}{n(n-1)} & \tfrac{b(b-1)}{n(n-1)}  & \dots  & \tfrac{b}{n}
    \end{bmatrix}, \notag
\end{eqnarray}
as claimed in \eqref{eq:prob_matrix_nice}. Therefore, 
\begin{eqnarray}\mP - pp^\top  = \begin{bmatrix}
    \tfrac{b}{n} -\frac{b^2}{n^2} & \tfrac{b(b-1)}{n(n-1)}& \dots  & \tfrac{b(b-1)}{n(n-1)} \\
    \tfrac{b(b-1)}{n(n-1)} & \tfrac{b}{n} & \dots  & \tfrac{b(b-1)}{n(n-1)} \\
    \vdots & \vdots & \ddots & \vdots \\
    \tfrac{b(b-1)}{n(n-1)} & \tfrac{b(b-1)}{n(n-1)}  & \dots  & \tfrac{b}{n}
    \end{bmatrix}, \notag
\end{eqnarray}

\item 
Letting $t=\frac{(a-1)k}{a(k-1)}$ and $s=1-t=\tfrac{k-a}{a(k-1)}$ the probability matrix of the approximate independent sampling satisfies
\begin{eqnarray}
\mP - pp^\top &=& \begin{bmatrix}
    p_{1}(1-p_1) & (t-1) p_1 p_2& \dots  & (t- 1) p_1 p_k  & 0 & \dots  & 0 \\
    (t- 1)  p_2 p_1 & p_{2}(1-p_2) & \dots  & (t- 1)  p_2 p_k & 0 & \dots  & 0\\
    \vdots & \vdots & \ddots & \vdots & 0 & \dots  & 0\\ 
   (t- 1)  p_n p_1 & (t- 1)  p_n p_2 & \dots  & p_{k}(1-p_k) & 0 & \dots  & 0 \\
   0 & 0 & \dots & 0 & 0 & \dots & 0 \\
   \vdots & \vdots & \vdots & \vdots & \vdots & \ddots & \vdots \\
  0& 0 & \dots  & 0 & 0 & \dots  & 0 
    \end{bmatrix} \notag \\
    &=& \Diag{p_1(1-p_1(1-s)),\dots, p_k(1-p_k(1-s)),0, \dots, 0} - s p_kp_k^\top \notag \\
    &\preceq& \Diag{p_1(1-p_1(1-s)),\dots, p_n(1-p_n(1-s)), 0, \dots, 0}, \notag
\end{eqnarray}
where $p_k = (p_1, \dots, p_k, 0, \dots, 0)^\top.$ 
Therefore, $v_i = 1-p_i(1-s)$ for $i \leq k$ and $v_i = 0$ otherwise works.

\item 
Finally, as remarked in the introduction, the standard uniform minibatch sampling ($b$--nice sampling) arises as a special case of the approximate independent sampling for the choice $p_i=b/n$. Thus $k = n,$  $a=b$ and hence $s=\tfrac{n-b}{b(n-1)}$. Based on the previous result, $v_i = 1 - \tfrac{b}{n}(1-\tfrac{n-b}{b(n-1)}) = \frac{n-b}{n-1}$ works.
\end{itemize}
\end{proof}

\section{Proof of Theorem~\ref{optimality}}

We first establish  a lemma we will need in order to prove Theorem~\ref{optimality}.

\begin{lemma} 
\label{opt_sampling}
Let  $0<L_1 \leq L_2\leq \dots \leq L_n$  be positive real numbers, $0<b\leq n$, and consider the optimization problem
\begin{eqnarray}
&\text{minimize}_{p \in \R^n} & \Omega(p) \eqdef \sum_{i =1}^n \frac{L_i^2}{p_i} \notag \\
&\text{subject to}&  \sum_{i =1}^n p_i = b, \label{eq:OPT_aux}\\
&& 0 \leq p_i \leq 1,\quad i=1,2,\dots,n. \notag 
\end{eqnarray}

Let  be the largest integer  for which 
$0<b +  k -n \leq \frac{\sum_{i =1}^k L_i}{L_k}$ (note that the inequality holds for $k=n-b+1$). Then \eqref{eq:OPT_aux} has the following solution:
\begin{equation}
\label{unique_sol}
p_i = 
\begin{cases}
       (b +  k - n)\frac{L_i}{\sum_{j=1}^k L_j}, &\quad\text{if } i \leq k,  \\
       1, &\quad\text{if } i > k.
\end{cases}
\end{equation} 
\end{lemma}

\begin{proof}
The Lagrangian of the problem is
\begin{eqnarray}
L(p,y,\lambda_1, ..,\lambda_n, u_1,...,u_n) = \sum_{i =1}^n \frac{L_i^2}{p_i} - \sum_{i =1}^n \lambda_ip_i - \sum_{i =1}^n u_i(1-p_i) + y\left(\sum_{i =1}^n p_i - b\right). \notag 
\end{eqnarray}
Th constraints are linear and hence  KKT conditions hold. The result can be deduced from the KKT conditions.
\end{proof}

We can now proceed with the proof. Since $n,b$ and $\bar{L}$ are  constants, the problem is equivalent to 
\begin{eqnarray*}
\text{minimize}_{\Sam} && \psi(\Sam)\eqdef \sum_{i = 1}^n \frac{v_iL_i^2}{p_i} \\
\text{subject to} && v_i \quad \text{satisfies} \quad \eqref{eq:ESO}.
\end{eqnarray*}
In view of \eqref{v_i}, 
\[\psi(\Sam) \overset{\eqref{v_i}}{\geq}  \sum_{i = 1}^n \frac{(1-p_i)L_i^2}{p_i} =  \sum_{i = 1}^n \frac{L_i^2}{p_i} -  \sum_{i = 1}^n L_i^2 = \Omega(p) -  \sum_{i = 1}^n L_i^2, \]
where function $\Omega(p)$ was defined in Lemma~\ref{opt_sampling}. Since $b = \Exp{|S|} = \sum_i p_i$, and $0\leq p_i \leq 1$ for all $i$, then in view of Lemma~\ref{opt_sampling} we have 
\[\Psi(\Sam) \geq \Omega(p^*) - \sum_{i = 1}^n L_i^2, \]
where $p^*$ is defined by \eqref{prob}.

On the other hand, from Lemma~\ref{lem:upperv} we know that the independent sampling $\Sam=\Sam^*$ with probability vector $p^*$ defined in \eqref{prob} satisfies inequality \eqref{eq:ESO} with  $v_i = 1-p_i$, and hence
\[\Psi(\Sam^*) =\Omega(p^*) -  \sum_{i = 1}^n L_i^2.\]
Hence, it is optimal. \qed

\section{Improvements}
\label{sec:improvements}
Let us compute $\alpha$ for uniform sampling.
\begin{eqnarray}
\alpha  &\overset{\eqref{def_K}}{=}&   \left(\frac{b}{n^2} \sum_{i=1}^n \frac{v_i L_i^2}{p_i}\right)/\bar{L}^2 \notag \\
	&\overset{\text{Lemma}~\ref{lem:upperv}}{=}&  \left(\frac{(n-b)}{(n-1)n} \sum_{i=1}^n  L_i^2\right)/\bar{L}^2 \notag \\
	&=& n\frac{(n-b)}{(n-1)}  \sum_{i=1}^n  L_i^2/ \left(\sum_{i=1}^n  L_i\right)^2 \notag 
\end{eqnarray}
It is easy to see that $L_{\max} \geq \bar{L}$. To prove that we have improved current best known rates, we need to show that $\bar{L}\alpha \leq L_{\max}$  and $\bar{L}^2\alpha \leq \frac{(n-b)}{(n-1)}L_{\max}^2$
\begin{proof}
\begin{eqnarray}
\bar{L}\alpha &=& n \frac{(n-b)}{(n-1)}\frac{\sum_{i=1}^n L_i^2}{\left( \sum_{i=1}^n L_i\right)^2} \bar{L} \leq \frac{\sum_{i=1}^n L_i^2}{\left( \sum_{i=1}^n L_i\right)}  =  \frac{\sum_{i=1}^n L_{\max}L_i}{\left( \sum_{i=1}^n L_i\right)} = L_{\max},  \notag \\
 \bar{L}^2\alpha &=& n \frac{(n-b)}{(n-1)}\frac{\sum_{i=1}^n L_i^2}{\left( \sum_{i=1}^n L_i\right)^2} \bar{L}^2 = \frac{(n-b)}{(n-1)}\frac{1}{n} \sum_{i=1}^n L_i^2 \leq  \frac{(n-b)}{(n-1)}\frac{1}{n} \sum_{i=1}^n L_{\max}^2 = \frac{(n-b)}{(n-1)}L_{\max}^2,  \notag
\end{eqnarray}
\end{proof}

Let's take $b = 1$. If  $L_n\gg L_i, \forall i \in [n]$, then $L_{\max} \approx n\bar{L}$ and $\alpha_{\Sam^*}\leq 1$ and $\alpha_{\Sam^u} \approx n$, which essentially means, that we can have in theory speedup by factor of $n$.

\section{Stochastic gradients evaluation complexity}
\label{sec:Rates}

\subsection{SVRG}
For \texttt{SVRG}, each outer loop costs $n+mb$ evaluations of stochastic gradient. If we want to obtain $\epsilon$-solution, following must hold (Theorem~\ref{thm:nonconvex-minibatch})
\begin{equation}
\frac{\alpha \bar{L}n^{(2/3)}(f(x^0)-f(x^*))}{bMm\nu_2} \leq \epsilon  \notag
\end{equation}

Combining these two equations with definition from Theorem~\ref{thm:nonconvex-minibatch}, we get  total complexity in terms of stochastic gradients evaluation
\begin{equation}
\frac{\mu_2 \bar{L}n^{(2/3)}(f(x^0)-f(x^*))}{\epsilon\nu_2}(1+\frac{\alpha}{3\mu_2}) \notag
\end{equation} 

\subsection{SAGA}
For \texttt{SAGA}, each loop costs $d+b$ evaluations of stochastic gradient. If we want to obtain $\epsilon$-solution, following must hold (Theorem~\ref{thm:nonconvex-minibatch_SAGA})
\begin{equation}
\frac{\alpha \bar{L}n^{(2/3)}(f(x^0)-f(x^*))}{bT\nu_2} \leq \epsilon \notag
\end{equation}

Combining these two equations with definition from Theorem~\ref{thm:nonconvex-minibatch_SAGA}, we get  total complexity in terms of stochastic gradients evaluation
\begin{equation}
 n + \frac{ \bar{L}n^{(2/3)}(f(x^0)-f(x^*))}{\epsilon\nu_3}(1+\alpha), \notag	
\end{equation}
because of evaluation of full gradient on the start. 
\subsection{SARAH}
For \texttt{SARAH} with one outer loot, each inner  loop costs $2b$ evaluations of stochastic gradient. If we want to obtain $\epsilon$-solution, following must hold (Theorem~\ref{thm:nonconvex_01_mb})
\begin{equation}
\frac{ 2\bar{L}(f(x^0)-f(x^*)) \left( \sqrt{1 + \frac{4m\alpha}{b}}\right)}{m} \leq \epsilon\notag
\end{equation}

Solving this equation for $m$, we get 
\begin{equation}
m \leq \frac{16 \alpha\bar{L}^2 (f(x^0)-f(x^*))^2 + \sqrt{16^2 \alpha^2\bar{L}^4 (f(x^0)-f(x^*))^4 + 16 \epsilon^2\bar{L}^2(f(x^0)-f(x^*))^2b^2}}{2b\epsilon^2} \notag
\end{equation}
Combining thise equation with complexity off each inner loop we obtain total complexity in terms of stochastic gradients evaluation
\begin{equation}
\frac{16 \alpha\bar{L}^2 (f(x^0)-f(x^*))^2 + \sqrt{16^2 \alpha^2\bar{L}^4 (f(x^0)-f(x^*))^4 + 16 \epsilon^2\bar{L}^2(f(x^0)-f(x^*))^2b^2}}{2\epsilon^2}. \notag
\end{equation}

\clearpage
\section{Proofs for SVRG}

\begin{lemma}
\label{lem:nonconvex-svrg}
For $c_t, c_{t+1}, \beta > 0$, suppose we have 
\begin{equation}
c_{t} = c_{t+1}(1 + \eta\beta + 2\eta^2K ) +K\eta^2\bar{L}.\notag
\end{equation}

Let $\eta$, $\beta$ and $c_{t+1}$ be chosen such that $\Gamma_t > 0$ (in Theorem~\eqref{thm:nonconvex-inter_is}). The iterate $x^{s+1}_t$ in Algorithm~\ref{alg:svrg_imp} satisfy the bound:
  \begin{eqnarray}
    \Exp{\|\nabla f(x^{s+1}_{t})\|^2} \leq \frac{R^{s+1}_{t} - R^{s+1}_{t+1}}{\Gamma_t},\notag
  \end{eqnarray}
  where $R^{s+1}_{t} \eqdef \Exp{f(x^{s+1}_{t}) + c_{t} \|x^{s+1}_{t} - \tilde{x}^{s}\|^2}$ for $0 \leq s \leq S-1$.
\end{lemma}
\begin{proof}
Since $f_i$ is $L_i$-smooth we have
\begin{eqnarray}
&\Exp{f_i(x^{s+1}_{t+1})} \leq \Exp{f_i(x^{s+1}_{t}) + \langle \nabla f_i(x^{s+1}_t), x^{s+1}_{t+1} - x^{s+1}_t \rangle + \tfrac{L_i}{2} \| x^{s+1}_{t+1} - x^{s+1}_t \|^2}.\notag
\end{eqnarray}
Summing through all $i$ and dividing by $n$ we obtain
\begin{eqnarray}
&\Exp{f(x^{s+1}_{t+1})} \leq \Exp{ f(x^{s+1}_{t}) + \langle \nabla f(x^{s+1}_t), x^{s+1}_{t+1} - x^{s+1}_t \rangle + \tfrac{\bar{L}}{2} \| x^{s+1}_{t+1} - x^{s+1}_t \|^2}.\notag
\end{eqnarray}
Using the \texttt{SVRG} update in Algorithm~\ref{alg:svrg_imp} and its unbiasedness ($\Exp{i_t}{v^{s+1}_t} = \nabla f(x_t^{s+1})$), the right hand side above is further upper bounded by
\begin{eqnarray}
  \Exp{f(x^{s+1}_{t}) - \eta \|\nabla f(x^{s+1}_{t})\|^2 + \tfrac{\bar{L}\eta^2}{2} \|v^{s+1}_t \|^2}.
\label{eq:svrg-proof-eq1}
\end{eqnarray}
Consider now the Lyapunov function
\begin{equation}
R^{s+1}_{t} \eqdef \Exp{ f(x^{s+1}_{t}) + c_{t} \|x^{s+1}_{t} - \tilde{x}^s\|^2 }. \notag
\end{equation}
For bounding it we will require the following:
\begin{eqnarray}
\Exp{\|x^{s+1}_{t+1} - \tilde{x}^s\|^2} &=& \Exp{ \|x^{s+1}_{t+1} - x^{s+1}_t + x^{s+1}_t - \tilde{x}^s\|^2} \notag \\
&=& \Exp{ \|x^{s+1}_{t+1} - x^{s+1}_t\|^2 + \|x^{s+1}_t - \tilde{x}^s\|^2} \notag \\
&& \quad+ 2\langle x^{s+1}_{t+1} - x^{s+1}_t, x^{s+1}_t - \tilde{x}^s\rangle] \notag \\
&=& \Exp{\eta^2\|v^{s+1}_t\|^2 + \|x^{s+1}_t - \tilde{x}^s\|^2} \notag \\
&& \quad - 2\eta \Exp{ \langle \nabla f(x^{s+1}_t), x^{s+1}_t - \tilde{x}^s\rangle } \notag \\
&\overset{\eqref{CS},\eqref{Young}}{\leq}& \Exp{ \eta^2\|v^{s+1}_t\|^2 + \|x^{s+1}_t - \tilde{x}^s\|^2} \notag \\
&& \qquad + 2 \eta \Exp{\tfrac{1}{2\beta} \|\nabla f(x^{s+1}_t)\|^2 + \tfrac{1}{2}\beta \| x^{s+1}_t - \tilde{x}^s \|^2 }.
\label{eq:svrg-proof-eq2}
\end{eqnarray}
The second equality follows from the unbiasedness of the update of \texttt{SVRG}. Plugging Equation~\eqref{eq:svrg-proof-eq1} and Equation~\eqref{eq:svrg-proof-eq2} into $R^{s+1}_{t+1}$, we obtain the following bound:
\begin{eqnarray}
R^{s+1}_{t+1} &\leq &\Exp{f(x^{s+1}_{t}) - \eta \|\nabla f(x^{s+1}_{t})\|^2 + \tfrac{\bar{L}\eta^2}{2} \|v^{s+1}_t \|^2}  \notag \\
&& \quad + \Exp{ c_{t+1}\eta^2\|v^{s+1}_t\|^2 + c_{t+1}\|x^{s+1}_t - \tilde{x}^s\|^2} \nonumber \\
&&  \quad + 2 c_{t+1}\eta \Exp{\tfrac{1}{2\beta} \|\nabla f(x^{s+1}_t)\|^2 + \tfrac{1}{2}\beta \| x^{s+1}_t - \tilde{x}^s \|^2 } \nonumber\\
&\leq & \Exp{f(x^{s+1}_{t}) - \left(\eta - \tfrac{c_{t+1}\eta}{\beta}\right) \|\nabla f(x^{s+1}_{t})\|^2}  + \left(\tfrac{\bar{L}\eta^2}{2} + c_{t+1}\eta^2 \right)\Exp{\|v^{s+1}_t\|^2} \nonumber\\
&&  \quad + \left( c_{t+1} + c_{t+1}\eta\beta \right) \Exp{ \| x^{s+1}_t - \tilde{x}^s \|^2 }.
\label{eq:svrg-proof-eq3}
\end{eqnarray}
To further bound this quantity, we use Lemma~\ref{lem:nonconvex-variance-lemma} to bound $\Exp{\|v^{s+1}_{t}\|^2}$, so that upon substituting it in Equation~\eqref{eq:svrg-proof-eq3}, we see that
\begin{eqnarray}
R^{s+1}_{t+1} &\overset{\eqref{nonconvex-variance-lemma}}{\leq} &\Exp{f(x^{s+1}_{t})} - \left(\eta - \tfrac{c_{t+1}\eta}{\beta} - \eta^2\bar{L} - 2c_{t+1}\eta^2\right) \Exp{\|\nabla f(x^{s+1}_{t})\|^2} \nonumber\\
&& \qquad + \left[c_{t+1}\bigl(1 + \eta\beta + 2\eta^2K\bigr)+\eta^2K\bar{L}\right] 
  \Exp{\| x^{s+1}_t - \tilde{x}^s \|^2 } \nonumber \\
& \leq & R^{s+1}_{t} - \bigl(\eta - \tfrac{c_{t+1}\eta}{\beta} - \eta^2\bar{L} - 2c_{t+1}\eta^2\bigr) \Exp{\|\nabla f(x^{s+1}_{t})\|^2}.
\end{eqnarray}
The second inequality follows from the definition of $c_{t}$ and $R_t^{s+1}$, thus concluding the proof. 
\end{proof}


\subsubsection*{Proof of Lemma~\ref{lem:nonconvex-svrg} and Theorem~\ref{thm:nonconvex-inter_is}}
\begin{proof}
  Using  Lemma~\ref{lem:nonconvex-svrg} and telescoping the sum, we obtain
\begin{eqnarray}
  \sum_{t=0}^{m-1} \Exp{\|\nabla f(x^{s+1}_{t})\|^2} \leq \frac{R^{s+1}_{0} - R^{s+1}_{m}}{\gamma_n}.
\end{eqnarray}
This inequality in turn implies that
\begin{equation}
  \label{eq:descent-property}
  \sum_{t=0}^{m-1} \Exp{\|\nabla f(x^{s+1}_{t})\|^2} \leq \frac{\Exp{f(\tilde{x}^s) - f(\tilde{x}^{s+1})}}{\gamma_n},
\end{equation}
where we used that $R^{s+1}_{m} = \Exp{f(x^{s+1}_m)} = \Exp{f(\tilde{x}^{s+1})}$ (since $c_m = 0$), and that $R^{s+1}_{0} = \Exp{f(\tilde{x}^s)}$ (since $x^{s+1}_0 = \tilde{x}^s$). Now sum over all epochs to obtain 
\begin{equation}
  \frac{1}{T} \sum_{s=0}^{S-1}\sum_{t=0}^{m-1} \Exp{\|\nabla f(x^{s+1}_{t})\|^2} \leq \frac{f(x^{0}) - f(x^*)}{T\gamma_n}.
\label{eq:nonconvex-cor-eq1}
\end{equation}
The above inequality used the fact that $\tilde{x}^0 = x^0$. Using the above inequality and the definition of $x_a$ in Algorithm~\ref{alg:svrg_imp}, we obtain the desired result.
\end{proof}

\subsubsection*{Proof of Theorem~\ref{thm:nonconvex-gen_is}}
\begin{proof}
  For our analysis, we will require an upper bound on $c_{0}$. Let $m = \lfloor Kn/(3\bar{L}^2 \rfloor \mu_0), \eta = \mu_0 \bar{L}/(Kn^{2/3}).$  We observe that $c_0 = \tfrac{\mu_0^2\bar{L}^3}{K n^{4/3}} \tfrac{(1 + \theta)^m - 1}{\theta}$ where $\theta = 2K\eta^2 + \eta\beta$. This is obtained using the relation $c_{t} = c_{t+1}(1 + \eta\beta + 2K\eta^2 ) +  \eta^2K\bar{L}$ and the fact that $c_m = 0$. Using the specified values of $\beta$ and $\eta$ we have
  \begin{eqnarray}
    \theta = 2K\eta^2+ \eta\beta = \frac{2\mu_0^2 \bar{L}^2}{K n^{4/3}} + \frac{\mu_0 \bar{L}^2}{K n} \leq \frac{3\mu_0\bar{L}^2}{K n}.
  \end{eqnarray}
  The above inequality follows since $\mu_0 \leq 1$ and $n \geq 1$. Using the above bound on $\theta$, we get 
  \begin{eqnarray}
    \label{eq:c0-bound}
    c_0 &=& \frac{\mu_0^2\bar{L}^3}{n^{2}K} \frac{(1 + \theta)^m - 1}{\theta} = \frac{\mu_0\bar{L} ((1 + \theta)^m - 1)}{2\mu_0 + n^{\tfrac13}} \nonumber \\
        &\leq & \frac{\mu_0\bar{L} ((1 + \frac{3\mu_0\bar{L}^2}{nK})^{\lfloor \nicefrac{Kn}{3\mu_0\bar{L}^2} \rfloor} - 1)}{2\mu_0 + n^{\tfrac13}} \nonumber \\
        &\leq & n^{-\tfrac13}(\mu_0\bar{L} (e - 1)),
  \end{eqnarray}
  wherein the second inequality follows upon noting that $(1 + \frac{1}{l})^l$ is increasing for $l>0$ and $\lim_{l \rightarrow \infty} (1 + \frac{1}{l})^l = e$ (here $e$ is the Euler's number). Now we can lower bound $\gamma_n$, as
  \begin{eqnarray}
    \gamma_n &= \min_t \bigl(\eta - \tfrac{c_{t+1}\eta}{\beta} - \eta^2\bar{L} - 2c_{t+1}\eta^2\bigr) \\ \notag
             &\geq \bigl(\eta - \tfrac{c_{0}\eta}{\beta} - \eta^2\bar{L} - 2c_{0}\eta^2\bigr) \geq \frac{\nu\bar{L}}{Kn^{\tfrac23}},\notag
  \end{eqnarray}
  where $\nu$ is a constant independent of $n$. The first inequality holds since $c_t$ decreases with $t$. The second inequality holds since (a) $c_0/\beta$ is upper bounded by a constant independent of $n$ as $c_0/\beta \leq \mu_0(e-1)$ (follows from Equation~\eqref{eq:c0-bound}), (b) $\eta^2\bar{L} \leq \mu_0 \eta$ and (c) $2c_0\eta^2 \leq 2\mu_0^2(e-1)\eta$ (follows from Equation~\eqref{eq:c0-bound}). By choosing $\mu_0$ (independent of $n$) appropriately, one can ensure that $\gamma_n \geq \nu \bar{L}/(Kn^{\tfrac23})$ for some universal constant $\nu$. For example, choosing $\mu_0 = 1/4$, we have $\gamma_n \geq \nu \bar{L}/(Kn^{\tfrac23})$ with $\nu = 1/40$. Substituting the above lower bound in Equation~\eqref{eq:nonconvex-cor-eq1}, we obtain the desired result. 
\end{proof}

\section{Minibatch SVRG}

\subsubsection*{Proof of Theorem~\ref{thm:nonconvex-minibatch}}

The proofs essentially follow along the lines of Lemma~\ref{lem:nonconvex-svrg}, Theorem~\ref{thm:nonconvex-inter_is} and Theorem~\ref{thm:nonconvex-gen_is} with the added complexity of mini-batch. We first prove few intermediate results before proceeding to the proof of Theorem~\ref{thm:nonconvex-minibatch}.

\begin{lemma}
  \label{lem:nonconvex-minibatch-svrg}
  Suppose we have
    \begin{eqnarray}
    \label{n312412asd}
     &\overline{R}^{s+1}_{t} \eqdef \Exp{f(x^{s+1}_{t}) + \overline{c}_{t} \|x^{s+1}_{t} - \tilde{x}^{s}\|^2},  \\\ \notag
     &\overline{c}_{t} = \overline{c}_{t+1}(1 + \eta\beta + \frac{2K\eta^2}{b} ) +  \frac{K\eta^2\bar{L}}{b},
     \end{eqnarray}
     for $\ 0 \leq s \leq S-1$ and  $0 \leq t \leq m-1$ and the parameters $\eta, \beta$ and $\overline{c}_{t+1}$ are chosen such that
\begin{equation}
 \left(\eta - \frac{\overline{c}_{t+1}\eta}{\beta} - \eta^2\bar{L} - 2\overline{c}_{t+1}\eta^2\right) \geq 0. 
 \notag
\end{equation}

     Then the iterates $x^{s+1}_t$ in the mini-batch version of Algorithm \ref{alg:svrg_imp} i.e., Algorithm \ref{alg:minibatch-svrg} with expected mini-batch size $b$ satisfy the bound:
  \begin{eqnarray}
    \Exp{\|\nabla f(x^{s+1}_{t})\|^2} \leq \frac{\overline{R}^{s+1}_{t} - \overline{R}^{s+1}_{t+1}}{\left(\eta - \frac{\overline{c}_{t+1}\eta}{\beta} - \eta^2\bar{L} - 2\overline{c}_{t+1}\eta^2\right)}, \notag
  \end{eqnarray}
\end{lemma}
\begin{proof}
Using essentially the same argument as the proof of Lemma \ref{lem:nonconvex-svrg} until Equation \eqref{eq:svrg-proof-eq3}, we have
\begin{eqnarray}
 \overline{R}^{s+1}_{t+1} &\leq &\Exp{(x^{s+1}_{t})}- \left(\eta - \tfrac{\overline{c}_{t+1}\eta}{\beta}\right) \|\nabla f(x^{s+1}_{t})\|^2 + \left(\tfrac{\bar{L}\eta^2}{2} + \overline{c}_{t+1}\eta^2 \right)\Exp{\|v^{s+1}_t\|^2 } \nonumber\\
&&   \quad + \left( \overline{c}_{t+1} + \overline{c}_{t+1}\eta\beta \right) \Exp{ \| x^{s+1}_t - \tilde{x}^s \|^2 }.
\label{eq:minibatch-svrg-proof-eq3}
\end{eqnarray}
We use Lemma~\ref{lem:nonconvex-minibatch-variance-lemma} in order to bound $\Exp{\|v^{s+1}_{t}\|^2}$ in the above inequality. Substituting it in Equation \eqref{eq:minibatch-svrg-proof-eq3}, we see that
\begin{eqnarray}
\overline{R}^{s+1}_{t+1} &\overset{\eqref{nonconvex-variance-lemma_b}}{\leq} & \Exp{f(x^{s+1}_{t})}  - \left(\eta - \tfrac{\overline{c}_{t+1}\eta}{\beta} - \eta^2\bar{L} - 2\overline{c}_{t+1}\eta^2\right) \Exp{\|\nabla f(x^{s+1}_{t})\|^2} \nonumber\\
&& \quad + \left[\overline{c}_{t+1}\bigl(1 + \eta\beta + \tfrac{2K\eta^2}{b}\bigr)+\tfrac{K\eta^2\bar{L}}{b}\right] 
  \Exp{\| x^{s+1}_t - \tilde{x}^s \|^2 } \nonumber \\
&\overset{\eqref{n312412asd}}{\leq} & \overline{R}^{s+1}_{t} - \bigl(\eta - \tfrac{\overline{c}_{t+1}\eta}{\beta} - \eta^2\bar{L} - 2\overline{c}_{t+1}\eta^2\bigr) \Exp{\|\nabla f(x^{s+1}_{t})\|^2}.
\end{eqnarray}
The second inequality follows from the definition of $\overline{c}_{t}$ and $\overline{R}^{s+1}_{t}$, thus concluding the proof. 
\end{proof}

The following theorem provides convergence rate of mini-batch\texttt{SVRG}.
\begin{theorem}
  Let $\overline{\gamma}_n$ denote the following quantity:
  \begin{eqnarray}
    \overline{\gamma}_n \eqdef \min_{0 \leq t \leq m-1}\quad\bigl(\eta - \tfrac{\overline{c}_{t+1}\eta}{\beta} - \eta^2\bar{L} - 2\overline{c}_{t+1}\eta^2\bigr). \notag
  \end{eqnarray}
  Suppose  $\overline{c}_m = 0$, $\overline{c}_{t} = \overline{c}_{t+1}(1 + \eta\beta + \frac{2K\eta^2}{b} ) +  \frac{K\eta^2\bar{L}}{b}$ for $t \in \{0, \dots, m-1\}$  and $\overline{\gamma}_n > 0$.  Then for the output $x_a$ of mini-batch version of Algorithm \ref{alg:svrg_imp} with mini-batch size $b$, we have
  \begin{eqnarray}
    \Exp{\|\nabla f(x_a)\|^2} \leq  \frac{f(x^{0}) - f(x^*)}{T\overline{\gamma}_n}, \notag
  \end{eqnarray}
  where $x^*$ is an optimal solution to \eqref{eq:1}.
  \label{thm:minibatch-nonconvex-inter}
\end{theorem}
\begin{proof}
  Using  Lemma \ref{lem:nonconvex-minibatch-svrg} and telescoping the sum, we obtain
\begin{eqnarray}
  \sum_{t=0}^{m-1} \Exp{\|\nabla f(x^{s+1}_{t})\|^2} \leq \frac{\overline{R}^{s+1}_{0} - \overline{R}^{s+1}_{m}}{\overline{\gamma}_n}. \notag
\end{eqnarray}
This inequality in turn implies that
\begin{eqnarray}
  \sum_{t=0}^{m-1} \Exp{\|\nabla f(x^{s+1}_{t})\|^2} \leq \frac{\Exp{f(\tilde{x}^s) - f(\tilde{x}^{s+1})}}{\overline{\gamma}_n}, \notag
\end{eqnarray}
where we used that $\overline{R}^{s+1}_{m} = \Exp{f(x^{s+1}_m)} = \Exp{f(\tilde{x}^{s+1})}$ (since $\overline{c}_m = 0$), and that $\overline{R}^{s+1}_{0} = \Exp{f(\tilde{x}^s)}$. Now sum over all epochs and using the fact that $\tilde{x}^0 = x^0$, we get the desired result.
\end{proof}

We now present the proof of Theorem \ref{thm:nonconvex-minibatch} using the above results.
\begin{proof}[Proof of Theorem \ref{thm:nonconvex-minibatch}]
  We first observe that using the specified values of $\beta = \bar{L}/ n^{1/3}$, $\eta = \mu_2b \bar{L}/(Kn^{2/3})$ and $\eta = \lfloor nK/(b\bar{L}^2\mu_2) \rfloor$we obtain
  \begin{eqnarray}
    \overline{\theta} \eqdef \frac{2K\eta^2}{b} + \eta\beta = \frac{2\mu_2^2b \bar{L}^2}{Kn^{4/3}} + \frac{\bar{L}^2\mu_2 b}{K n} \leq \frac{3\mu_2\bar{L}^2 b}{K n}. \notag
  \end{eqnarray}
  The above inequality follows since $\mu_2 \leq 1$ and $n \geq 1$. For our analysis, we will require the following bound on $\overline{c}_{0}$:
  \begin{eqnarray}
    \label{eq:minibatch-c0-bound}
    \overline{c}_0 &=& \frac{\mu_2^2b^2\bar{L}^3}{K bn^{4/3}} \frac{(1 + \overline{\theta})^m - 1}{\overline{\theta}} = \frac{\mu_2b\bar{L} ((1 + \overline{\theta})^m - 1)}{2b\mu_2 +  bn^{1/3}} \nonumber \\
        &\leq& n^{-1/3}(\mu_2\bar{L} (e - 1)),
  \end{eqnarray}
  wherein the first equality holds due to the relation $\overline{c}_{t} = \overline{c}_{t+1}(1 + \eta\beta + \tfrac{2K\eta^2}{b}) +  \tfrac{K\eta^2\bar{L}}{b}$, and the inequality follows upon again noting that $(1 + 1/l)^l$ is increasing for $l>0$ and  $\lim_{l \rightarrow \infty} (1 + \frac{1}{l})^l = e$. Now we can lower bound $\overline{\gamma}_n$, as
  \begin{eqnarray}
    \overline{\gamma}_n &=& \min_t \bigl(\eta - \tfrac{\overline{c}_{t+1}\eta}{\beta} - \eta^2\bar{L} - 2\overline{c}_{t+1}\eta^2\bigr)  \notag \\ 
             &\geq &\bigl(\eta - \tfrac{\overline{c}_{0}\eta}{\beta} - \eta^2\bar{L} - 2\overline{c}_{0}\eta^2\bigr) \geq \frac{b\bar{L}\nu_2}{Kn^{2/3}}, \notag
  \end{eqnarray}
  where $\nu_2$ is a constant independent of $n$. The first inequality holds since $\overline{c}_t$ decreases with $t$. The second one holds since (a) $\overline{c}_0/\beta$ is upper bounded by a constant independent of $n$ as $\overline{c}_0/\beta \leq \mu_2(e-1)$ (due to Equation\eqref{eq:minibatch-c0-bound}), (b) $\eta^2\bar{L} \leq \mu_2\eta$ (as $b \leq \nicefrac{K}{\bar{L}^2n^{2/3}}$) and (c) $2\overline{c}_{0}\eta^2 \leq 2\mu_2^2(e-1)\eta$ (again due to Equation \eqref{eq:minibatch-c0-bound} and the fact $b \leq \nicefrac{K}{\bar{L}^2n^{2/3}}$). By choosing an appropriately small constant $\mu_2$ (independent of n), one can ensure that $\overline{\gamma}_n \geq {\bar{L}b\nu_2}/(Kn^{2/3})$ for some universal constant $\nu_2$. For example, choosing $\mu_2 = 1/4$, we have $\overline{\gamma}_n \geq {\bar{L}b\nu_2}/(Kn^{2/3})$ with $\nu_2 = 1/40$. Substituting the above lower bound in Theorem \ref{thm:minibatch-nonconvex-inter}, we obtain the desired result. 
\end{proof}

\subsubsection*{Lemmas}
\begin{lemma}
\label{lem:nonconvex-variance-lemma}
  For the intermediate iterates $v^{s+1}_t$ computed by Algorithm~\ref{alg:svrg_imp}, we have the following:
 \begin{eqnarray}
  \label{nonconvex-variance-lemma}
    \Exp{ \|v^{s+1}_t\|^2 } \leq 2\Exp{\|\nabla f(x^{s+1}_{t})\|^2} + 2K \Exp{\|x^{s+1}_{t} - \tilde{x}^{s}\|^2}. 
\end{eqnarray}
\end{lemma}
\begin{proof}
The proof simply follows from the proof of Lemma~\ref{lem:nonconvex-minibatch-variance-lemma} with $\Sam_t = \{i_t\}$.
\end{proof}

We now present a result to bound the variance of mini-batch \texttt{SVRG}.

\begin{lemma}
\label{lem:nonconvex-minibatch-variance-lemma}
  Let $v^{s+1}_t$ be computed by the mini-batch version of Algorithm~\ref{alg:svrg_imp} i.e., Algorithm~\ref{alg:minibatch-svrg} with sampling $\Sam$. Then,
  \begin{eqnarray}
  \label{nonconvex-variance-lemma_b}
    \Exp{\|v^{s+1}_t\|^2} \leq 2\Exp{\|\nabla f(x^{s+1}_{t})\|^2} + \tfrac{2K}{b} \Exp{\|x^{s+1}_{t} - \tilde{x}^{s}\|^2}.
\end{eqnarray}
\end{lemma}
\begin{proof}
For the simplification, we use the following notation:
\begin{equation}
\zeta_{t}^{s+1} = \sum_{i_t \in \Sam_t} \frac{1}{np_{i_t}}\left(\nabla f_{i_t}(x^{s+1}_{t}) - \nabla f_{i_t}(\tilde{x}^{s})\right). \notag
\end{equation}

We use the definition of $v^{s+1}_t$ to get
\begin{eqnarray}
\Exp{\|v^{s+1}_t\|^2} &=& \Exp{\|\zeta_{t}^{s+1} + \nabla f(\tilde{x}^{s}) \|^2 } \notag \\
&=& \Exp{\| \zeta_{t}^{s+1} + \nabla f(\tilde{x}^{s}) - \nabla f(x^{s+1}_{t}) +  \nabla f(x^{s+1}_{t}) \|^2} \notag \\
&\leq & 2\Exp{\|\nabla f(x^{s+1}_{t})\|^2} + 2 \Exp{\|\zeta_{t}^{s+1} - \Exp{\zeta_{t}^{s+1}} \|^2} \notag \\
&= & 2\Exp{\|\nabla f(x^{s+1}_{t})\|^2} \notag \\
&& \quad+2\Exp{\left\|\sum_{i_t \in \Sam_t} \left(\frac{1}{np_{i_t}}\left(\nabla f_{i_t}(x^{s+1}_{t}) - \nabla f_{i_t}(\tilde{x}^{s})\right) - \Exp{\zeta_{t}^{s+1}} \right)  \right\|^2 }.\notag
\end{eqnarray}
The first inequality follows from fact that $\norm{x+y}^2 \leq 2 \norm{x}^2 +2\norm{y}^2$ and the fact that $ \Exp{\zeta_{t}^{s+1}} = \nabla f(x^{s+1}_{t}) - \nabla f(\tilde{x}^{s})$. From the above inequality, we get
\begin{eqnarray}
\Exp{\|v^{s+1}_t\|^2}
& \overset{\eqref{lem:upperv}}{\leq} & 2\Exp{\|\nabla f(x^{s+1}_{t})\|^2} + 2 \sum_{i=1}^n\frac{v_i p_i}{n^2p_i^2} \left\|\left(\nabla f_{i}(x^{s+1}_{t}) - \nabla f_{i}(\tilde{x}^{s})\right)\right\|^2 \notag \\
& \overset{\eqref{L_i-smooth}, \eqref{def_K}}{\leq} & 2\Exp{\|\nabla f(x^{s+1}_{t})\|^2} + \frac{2K}{b} \Exp{\|x^{s+1}_{t} - \tilde{x}^{s}\|^2}.\notag
\end{eqnarray}
\end{proof}

\clearpage
\section{Proofs for SAGA}

\begin{lemma}
\label{lem:nonconvex-saga}
For $c_t, c_{t+1}, \beta > 0$, suppose we have 
\begin{equation}
c_{t} = c_{t+1}(1 - \tfrac{d}{ n} +  \eta\beta + 2\frac{K\eta^2}{b} ) + \frac{K \eta^2\bar{L}}{b}. \notag
\end{equation}

Also let $\eta$, $\beta$ and $c_{t+1}$ be chosen such that $\Gamma_{t} > 0$. Then, the iterates $\{x^t\}$ of Algorithm~\ref{alg:saga} satisfy the bound
\begin{eqnarray}
\Exp{\|\nabla f(x^{t})\|^2} \leq \frac{R^{t} - R^{t+1}}{\Gamma_t}, \notag
\end{eqnarray}
where $R^{t} \eqdef \Exp{f(x^{t})} +  c_{t} \max_{i \in [n]} \Exp{\|x^{t} - \alpha_i^{t}\|^2}$.
\end{lemma}

\begin{proof}
Since $f$ is $\bar{L}$-smooth we have
\begin{eqnarray}
&\Exp{f(x^{t+1})} \leq \Exp{ f(x^{t}) + \langle \nabla f(x^t), x^{t+1} - x^t \rangle + \tfrac{\bar{L}}{2} \| x^{t+1} - x^t \|^2}. \notag
\end{eqnarray}
We first note that the update in Algorithm~\ref{alg:saga} is unbiased i.e., $\Exp{v^t} = \nabla f(x^t)$. By using this property of the update on the right hand side of the inequality above, we get the following:
\begin{equation}
\Exp{f(x^{t+1})} \leq  \Exp{f(x^{t}) - \eta \|\nabla f(x^{t})\|^2 + \tfrac{\bar{L}\eta^2}{2} \|v^t \|^2}.
\label{eq:saga-proof-eq1}
\end{equation}
Here we used the fact that $x^{t+1} - x^{t} = -\eta v^t$ (see Algorithm~\ref{alg:minibatch-saga}). Consider now the Lyapunov function
\begin{equation}
R^{t} \eqdef \Exp{f(x^{t})} +  c_{t} \max_{i \in [n]} \Exp{\|x^{t} - \alpha_i^{t}\|^2}. \notag
\end{equation}
For bounding $R^{t+1}$ we need the following:
\begin{equation}
\label{eq:aux-term}
 \Exp{ \|x^{t+1} - \alpha_{i}^{t+1}\|^2} =   \frac{d}{ n} \Exp{\|x^{t+1} - x^{t}\|^2} +  \frac{n-d}{n} \underbrace{\Exp{\|x^{t+1} - \alpha_{i}^{t}\|^2}}_{T_1},
\end{equation}

The above equality follows from the definition of $\alpha^{t+1}_i$ and the definition of randomness of index $j_t$ in Algorithm~\ref{alg:saga} and Algorithm~\ref{alg:minibatch-saga}. The term $T_1$ in~\eqref{eq:aux-term} can be bounded as follows
\begin{eqnarray}
T_1 &=& \Exp{\|x^{t+1} - x^t + x^t - \alpha_{i}^{t}\|^2} \nonumber \\
&=& \Exp{ \|x^{t+1} - x^t\|^2 + \|x^t - \alpha_{i}^{t}\|^2 } + 2\langle x^{t+1} - x^t, x^t - \alpha_{i}^{t}\rangle] \nonumber \\
&=& \Exp{\|x^{t+1} - x^t\|^2 + \|x^t - \alpha_{i}^{t}\|^2} - 2\eta \Exp{\langle \nabla f(x^t), x^t -\alpha_{i}^{t}\rangle } \nonumber \\
& \overset{\eqref{CS},\eqref{Young}}{\leq} & \Exp{\|x^{t+1} - x^t\|^2 + \|x^t - \alpha_{i}^{t}\|^2}  + 2 \eta \Exp{\tfrac{1}{2\beta} \|\nabla f(x^t)\|^2 + \tfrac{1}{2}\beta \| x^t - \alpha_{i}^{t}\|^2 }\nonumber \\
& \leq & \Exp{\|x^{t+1} - x^t\|^2} +\max_{i \in [n]} \Exp{ \|x^t - \alpha_{i}^{t}\|^2}  + 2 \eta \Exp{\tfrac{1}{2\beta} \|\nabla f(x^t)\|^2} + \eta\beta \max_{i \in [n]} \Exp{\| x^t - \alpha_{i}^{t}\|^2 }.
\label{eq:saga-proof-eq2}
\end{eqnarray}
The second equality again follows from the unbiasedness of the update of \texttt{SAGA}. The last inequality follows from a simple application of Cauchy-Schwarz and Young's inequality. Plugging~\eqref{eq:saga-proof-eq1} and~\eqref{eq:saga-proof-eq2} into $R^{t+1}$, we obtain the following bound:
\begin{eqnarray}
 R^{t+1} &\leq & \Exp{f(x^{t}) - \eta \|\nabla f(x^{t})\|^2 + \tfrac{\bar{L}\eta^2}{2} \|v^t \|^2}  \nonumber \\
 && \quad + \Exp{c_{t+1}\|x^{t+1} - x^t\|^2} + c_{t+1}\frac{n-d}{n}  \max_{i \in [n]} \Exp{\|x^t - \alpha_i^t\|^2} \nonumber \\
&& \quad  + \frac{2(n-1)c_{t+1}\eta}{n}  \Exp{\tfrac{1}{2\beta} \|\nabla f(x^t)\|^2} + \tfrac{1}{2}\beta \max_{i \in [n]} \Exp{ \| x^t - \alpha_i^t \|^2 } \nonumber\\
&\leq &\Exp{ f(x^{t}) - \left(\eta - \tfrac{c_{t+1}\eta}{\beta}\right) \|\nabla f(x^{t})\|^2} + \left(\tfrac{\bar{L}\eta^2}{2} + c_{t+1}\eta^2 \right)\Exp{\|v^t\|^2} \nonumber\\
&&  \quad + \left(\frac{n-d}{n} c_{t+1} + c_{t+1}\eta\beta \right)  \max_{i \in [n]} \Exp{\| x^t - \alpha_i^t\|^2 },
\label{eq:saga-proof-eq3}
\end{eqnarray}
where we use that  $\| x^t - \alpha_{i_{\max}}^{t}\|^2 \leq  \max_{i \in [n]} \|x^t - \alpha_i^t\|^2$
To further bound the quantity in~\eqref{eq:saga-proof-eq3}, we use Lemma~\ref{lem:nonconvex-variance-lemma_2} to bound $\Exp{\|v^{t}\|^2}$, so that upon substituting it into~\eqref{eq:saga-proof-eq3}, we obtain
\begin{eqnarray}
 R^{t+1} &\overset{\eqref{nonconvex-variance-lemma_2}}{\leq} & \Exp{f(x^{t})} - \left(\eta - \frac{c_{t+1}\eta}{\beta} - \eta^2\bar{L} - 2c_{t+1}\eta^2\right) \Exp{\|\nabla f(x^{t})\|^2} \nonumber\\
&& \quad + \left[c_{t+1}\bigl(1 - \frac{d}{ n} + \eta\beta + 2\frac{K\eta^2}{b}\bigr)+\frac{K\eta^2\bar{L}}{b}\right]
 \max_{i \in [n]} \Exp{  \| x^t - \alpha_i^t \|^2 } \nonumber \\
& \leq & R^{t} - \bigl(\eta - \tfrac{c_{t+1}\eta}{\beta} - \eta^2\bar{L} - 2c_{t+1}\eta^2\bigr) \Exp{ \|\nabla f(x^{t})\|^2 }.
\end{eqnarray}
The second inequality follows from the definition of $c_{t}$ i.e., $c_{t} = c_{t+1}\bigl(1 - \tfrac{d}{n} + \eta\beta + 2\frac{K\eta^2}{b}\bigr)+\frac{K\eta^2\bar{L}}{b}$ and $R^t$ specified in the statement, thus concluding the proof.
\end{proof}

The following lemma provides a bound on the variance of the update used in Minibatch \texttt{SAGA} algorithm. More specifically, it bounds the quantity $\Exp{\|v^t\|^2}$.

\begin{lemma}
\label{lem:nonconvex-variance-lemma_2}
Let $v^t$ be computed by Algorithm~\ref{alg:minibatch-saga}. Then,
\begin{eqnarray}
\label{nonconvex-variance-lemma_2}
\Exp{\|v^t\|^2} \leq 2\Exp{\|\nabla f(x^{t})\|^2} + \frac{2K}{b}  \max_{i \in [n]}\Exp{\|x^{t} - \alpha_i^t\|^2}.
\end{eqnarray}
\end{lemma}
\begin{proof}
For ease of exposition, we use the notation
\begin{equation}
\zeta_i^{t} \eqdef \frac{1}{np_i}\left(\nabla f_{i}(x^{t}) - \nabla f_{i}(\alpha_{i}^{t})\right). \notag
\end{equation}

Using the convexity of $\|\!\cdot\!\|^2$ and the definition of $v^t$ we get
\begin{eqnarray}
\Exp{\|v^t\|^2} &=& \Exp{\|\sum_{i \in \Sam_t}\zeta_i^{t} + \tfrac{1}{n}\sum_{i=1}^n\nabla f(\alpha_i^t) \|^2}  \notag \\
&=& \Exp{\| \sum_{i \in \Sam_t}\zeta_i^{t}  + \tfrac{1}{n}\sum_{i=1}^n\nabla f(\alpha_i^t) - \nabla f(x^{t}) +  \nabla f(x^{t}) \|^2} \notag \\
&\leq & 2\Exp{\|\nabla f(x^{t})\|^2} + 2 \Exp{\|\sum_{i \in \Sam_t}\zeta_i^{t} - \Exp{\zeta^{t}}\|^2} \notag \\
&\overset{\eqref{lem:upperv}}{\leq} & 2\Exp{\|\nabla f(x^{t})\|^2} + 2 \sum_{i=1}^n \Exp{p_{i_t} \|\zeta_{i_t}^{t}\|^2}.\notag
\end{eqnarray}
The first inequality follows from the fact that $\|a + b\|^2 \leq 2(\|a\|^2 + \|b\|^2)$ and that $ \Exp{\zeta^{t}} = \nabla f(x^{t}) - \tfrac{1}{n} \sum_{i=1}^n \nabla f(\alpha_i^{t})$.v
\begin{eqnarray}
\Exp{\|v^t\|^2} & \leq &2\Exp{\|\nabla f(x^{t})\|^2} + 2 \sum_{i=1}^n \Exp{\frac{p_{i}}{n^2 p_{i}^2} \|\nabla f_{i}(x^{t}) - \nabla f_{i}(\alpha_i^{t})\|^2 } \notag \\
& \overset{\eqref{L_i-smooth}, \eqref{eq:key_inequality}}{\leq} &2\Exp{ \|\nabla f(x^{t})\|^2} +  2\sum_{i=1}^n \Exp{ \frac{v_i L_{i}^2}{n^2 p_{i}} \|x^{t} - \alpha_i^{t}\|^2} \notag \\
& \overset{\eqref{def_K}}{\leq} &2\Exp{ \|\nabla f(x^{t})\|^2} + \frac{2K}{b}  \max_{i \in [n]}\Exp{ \|x^{t} - \alpha_i^{t}\|^2}.
\end{eqnarray}
The last inequality follows from $L_{i}$-smoothness of $f_{i}$ and using properties of $\Sam$ sampling, thus concluding the proof.
\end{proof}

\subsubsection*{Proof of Theorem~\ref{thm:nonconvex-inter}}
\begin{proof}
  We apply telescoping sums to the result of Lemma~\ref{lem:nonconvex-saga} to obtain
\begin{eqnarray}
   \gamma_n \sum_{t=0}^{T-1} \Exp{\|\nabla f(x^{t})\|^2} &\leq \sum_{t=0}^{T-1} \Gamma_t \Exp{\|\nabla f(x^{t})\|^2} \leq R^{0} - R^{T}.
\end{eqnarray}
The first inequality follows from the definition of $\gamma_n$. This inequality in turn implies the bound
\begin{equation}
  \label{eq:descent-property_2}
  \sum_{t=0}^{T-1} \Exp{\|\nabla f(x^{t})\|^2} \leq \frac{\Exp{f(x^0) - f(x^{T})}}{\gamma_n},
\end{equation}
where we used that $R^{T} = \Exp{f(x^{T})}$ (since $c_T = 0$), and that $R^{0} = \Exp{f(x^0)}$ (since $\alpha_i^0 = x^0$ for $i \in [n]$). Using inequality~\eqref{eq:descent-property_2}, the optimality of $x^*$, and the definition of $x_a$ in Algorithm~\ref{alg:saga}, we obtain the desired result.
\end{proof}

\subsubsection*{Proof of Theorem~\ref{thm:nonconvex-gen} and Theorem~\ref{thm:nonconvex-minibatch_SAGA} }
\begin{proof}
With the values of $\mu_3 = 1/3, \nu_3 = 12$ $\eta = b\bar{L}/(3Kn^{2/3})$, $d = b\bar{L}^2/K$ and $\beta = \bar{L}/n^{1/3}$, let us first establish an upper bound on $c_t$. Let $\theta$ denote $\tfrac{\bar{L}^2b}{K n} - \eta\beta - 2K\eta^2/b$. Observe that $\theta < 1$ and $\theta \geq 4\bar{L}^2b/(9K n)$. This is due to the specific values of $\eta$ and $\beta$ and lower bound of $K$. Also, we have $c_{t} = c_{t+1}(1 - \theta) +  K\eta^2\bar{L}/b$. Using this relationship, it is easy to see that $c_{t} = K\eta^2\bar{L} \tfrac{1 - (1 - \theta)^{T-t}}{b\theta}$. Therefore, we obtain the bound
\begin{equation}
\label{eq:c-t-upperbound}
c_{t}  =K \eta^2\bar{L} \frac{1 - (1 - \theta)^{T-t}}{b\theta} \leq \frac{K\eta^2\bar{L}}{b\theta} \leq \frac{\bar{L}}{4 n^{1/3}},
\end{equation}
for all $0 \leq t \leq T$, where the inequality follows from the definition of $\eta$ and the fact that $\theta \geq 4\bar{L}^2b/(9K n)$. Using the above upper bound on $c_t$ we can conclude that
\begin{eqnarray}
\gamma_n = \min_t \left(\eta - \frac{c_{t+1}\eta}{\beta} - \eta^2\bar{L} - 2c_{t+1}\eta^2\right) \geq \frac{\bar{L}b}{12Kn^{2/3}}, \notag
\end{eqnarray}
upon using the following inequalities: (i) $c_{t+1}\eta/\beta \leq \eta/4$, (ii) $\eta^2L \leq \eta/3$ and (iii) $2c_{t+1}\eta^2 \leq \eta/6$, which hold due to the upper bound on $c_t$ in~\eqref{eq:c-t-upperbound} and if $b \leq \nicefrac{K}{\bar{L}^2}n^{2/3}$. Substituting this bound on $\gamma_n$ in Theorem~\ref{thm:nonconvex-inter},  we obtain the desired result. 
\end{proof}

Theorem~\ref{thm:nonconvex-gen} is special case with $ b = 1$ and $d = 1$.

\subsection*{SARAH-non-convex}
 This lemmas are modification of lemmas appeared in \cite{SARAH-nonconvex} for importance sampling with mini-batch.

\begin{lemma}\label{lem_main_derivation_mb}
Consider \texttt{SARAH}, then we have 
\begin{eqnarray}
 \label{eq:001}
\sum_{t=0}^{m} \Exp{ \| \nabla f(x^{t})\|^2 } &\leq & \frac{2}{\eta} [ f(x^{0}) - f(x^{*})] + \sum_{t=0}^{m} \Exp{ \| \nabla f(x^{t}) - v^{t} \|^2 }  \notag \\
&& \quad - ( 1 - \bar{L}\eta ) \sum_{t=0}^{m} \Exp{  \| v^{t} \|^2 }, 
\end{eqnarray}
where $x_{*}$ is an optimal solution of \eqref{eq:1}. 
\end{lemma}

\begin{proof}
By $\bar{L}$-smoothness of $f$ and $x^{t+1} = x^{t} - \eta v^{t}$, we have
\begin{eqnarray}
\Exp{f(x^{t+1})} & \leq &  \Exp{ f(x^{t})}- \eta \Exp{\nabla f(x^{t})^\top v^{t}}
+ \frac{\bar{L}\eta^2}{2} \Exp{ \| v^{t} \|^2 } \notag
\\
& =& \Exp{f(x^{t})}- \frac{\eta}{2} \Exp{ \| \nabla f(x^{t})\|^2 }
+ \frac{\eta}{2} \Exp{\| \nabla f(x^{t}) - v^{t} \|^2 }\notag \\
&& \quad - \left( \frac{\eta}{2} - \frac{\bar{L}\eta^2}{2} \right) \Exp{ \| v^{t} \|^2 },\notag
\end{eqnarray}
where the last equality follows from the fact
$r^\top q = \frac{1}{2}\left[\|r\|^2 + \|q\|^2 - \|r-q\|^2\right],$ for any $r,q\in \R^d$.

By summing over $t = 0,\dots,m$, we have
\begin{eqnarray}
\Exp{f(x^{m + 1})}& \leq &  \Exp{f(x^{0})}- \frac{\eta}{2} \sum_{t=0}^{m} \Exp{ \| \nabla f(x^{t})\|^2 }+ \frac{\eta}{2} \sum_{t=0}^{m} \Exp{ \| \nabla f(x^{t}) - v^{t} \|^2 } \notag \\
&& \quad- \left( \frac{\eta}{2} - \frac{L\eta^2}{2} \right) \sum_{t=0}^{m} \Exp{  \| v^{t} \|^2 },  \notag
\end{eqnarray}
which is equivalent to ($\eta>0$):
\begin{eqnarray}
\sum_{t=0}^{m} \Exp{ \| \nabla f(x^{t})\|^2 }  & \leq& \frac{2}{\eta} \Exp{f(x^{0}) - f(x^{m+1})} + \sum_{t=0}^{m} \Exp{ \| \nabla f(x^{t}) - v^{t} \|^2 }   \notag\\
&& \quad - ( 1 - \bar{L}\eta ) \sum_{t=0}^{m} \Exp{ \| v^{t} \|^2 } \notag \\
& \leq& \frac{2}{\eta} [ f(x^{0}) - f(x^{*})] + \sum_{t=0}^{m} \Exp{ \| \nabla f(x^{t}) - v^{t} \|^2 }\notag \\ 
 && \quad- ( 1 - \bar{L}\eta ) \sum_{t=0}^{m} \Exp{ \| v^{t} \|^2},    \notag
\end{eqnarray}

where the last inequality follows since $x^{*}$ is an optimal solution of \eqref{eq:1}. (Note that $x^{0}$ is given.) 

\end{proof}

\begin{lemma}\label{lem:var_diff_mb}
 Consider $v^{t}$ defined in \texttt{SARAH}, then for any $t\geq 1$, 
\begin{eqnarray}
\Exp{\| \nabla f(x^{t}) - v^{t} \|^2 }
= \sum_{j = 1}^{t} \Exp{ \| v^{j} - v^{j-1} \|^2 } 
 - \sum_{j = 1}^{t} \Exp{ \| \nabla f(x^{j}) - \nabla f(x^{j-1}) \|^2 }. \notag
\end{eqnarray}
\end{lemma}

\begin{proof}
Let $\mathcal{F}_{j} = \sigma(x^0,i_1,i_2,\dots,i_{j-1})$ be the $\sigma$-algebra generated by $x^0,i_1,i_2,\dots,i_{j-1}$; $\mathcal{F}_{0} = \mathcal{F}_{1} = \sigma(x^0)$. Note that $\mathcal{F}_{j}$ also contains all the information of $x^0,\dots,x^{j}$ as well as $v^0,\dots,v^{j-1}$. For $j \geq 1$, we have
\begin{eqnarray}
\Exp{\| \nabla f(x^{j}) - v^{j} \|^2 | \mathcal{F}_{j} ]}
& = & \Exp{ \| [\nabla f(x^{j-1}) - v^{j-1} } + [ \nabla f(x^{j}) - \nabla f(x^{j-1}) ] \notag \\
&& \quad - [ v^{j} - v^{j-1} ] \|^2 | \mathcal{F}_{j} ] \notag
\\
& = &\| \nabla f(x^{j-1}) - v^{j-1} \|^2 + \| \nabla f(x^{j}) - \nabla f(x^{j-1}) \|^2 \notag \\ 
&& \quad+ \Exp{ \| v^{j} - v^{j-1}  \|^2 | \mathcal{F}_{j} } \notag
\\
&&\quad + 2 ( \nabla f(x^{j-1}) - v^{j-1} )^\top ( \nabla f(x^{j}) - \nabla f(x^{j-1}) ) \notag\\
&&\quad - 2 ( \nabla f(x^{j-1}) - v^{j-1} )^\top \Exp{v^{j} - v^{j-1} | \mathcal{F}_{j} }\notag  \\
&&\quad - 2 ( \nabla f(x^{j}) - \nabla f(x^{j-1}) )^\top \Exp{ v^{j} - v^{j-1} | \mathcal{F}_{j} } \notag
\\
& = &\| \nabla f(x^{j-1}) - v^{j-1} \|^2 - \| \nabla f(x^{j}) - \nabla f(x^{j-1}) \|^2 \notag\\
&& \quad + \Exp{ \| v^{j} - v^{j-1}  \|^2 | \mathcal{F}_{j} }, \notag
\end{eqnarray}
where the last equality follows from
\begin{eqnarray}
\Exp{ v^{j} - v^{j-1} | \mathcal{F}_{j} }& = &\Exp{ \sum_{i \in I_{j}} \frac{1}{np_i}\nabla f_{i} (x^{j}) - \nabla f_{i}(x^{j-1})] \Big | \mathcal{F}_{j} }\notag\\
&=&  \sum_{i=1}^{n} \frac{p_i}{np_i}[\nabla f_{i} (x^{j}) - \nabla f_{i}(x^{j-1})] = \nabla f(x^{j}) - \nabla f(x^{j-1}).\notag
\end{eqnarray}

By taking expectation for the above equation, we have
\begin{eqnarray}
\Exp{\| \nabla f(x^{j}) - v^{j} \|^2 }&=& \Exp{\| \nabla f(x^{j-1}) - v^{j-1} \|^2 } - \Exp{ \| \nabla f(x^{j}) - \nabla f(x^{j-1}) \|^2 } \notag \\
&& \quad+ \Exp{\| v^{j} - v^{j-1} \|^2 }. \notag
\end{eqnarray}

Note that $\| \nabla f(x^{0}) - v^{0} \|^2 = 0$. By summing over $j = 1,\dots,t\ (t\geq 1)$, we have
\begin{eqnarray}
\Exp{\| \nabla f(x^{t}) - v^{t} \|^2 }  = \sum_{j = 1}^{t} \Exp{\| v^{j} - v^{j-1} \|^2 } - \sum_{j = 1}^{t} \Exp{ \| \nabla f(x^{j}) - \nabla f(x^{j-1}) \|^2 }. \notag
\end{eqnarray}

\end{proof}

With the above Lemmas, we can derive the following upper bound for $\Exp{ \| \nabla f(x^{t}) - v^{t} \|^2 }$. 

\begin{lemma}\label{lem:var_diff_mb_02}
 Consider $v^{t}$ defined in \texttt{SARAH}. Then for any $t\geq 1$, 
\begin{eqnarray}
\Exp{ \| \nabla f(x^{t}) - v^{t} \|^2 }  \leq \frac{1}{b} K \eta^2 \sum_{j=1}^{t} \Exp{\| v^{j-1} \|^2}. \notag
\end{eqnarray}
\end{lemma}

\begin{proof}
Let
\begin{eqnarray}
\label{xi_value}
\xi_t = \frac{1}{np_i}\left(\nabla f_{t} (x^{j}) - \nabla f_{t}(x^{j-1})\right)
\end{eqnarray}

We have
\begin{eqnarray}
 &&\Exp{\| v^j - v^{j-1} \|^2 | \mathcal{F}_{j} } - \| \nabla f(x^{j}) - \nabla f(x^{j-1}) \|^2\notag \\
& =& \Exp{  \Big\|  \sum_{i \in I_{j}} \frac{1}{np_i}[\nabla f_{i} (x^{j}) - \nabla f_{i}(x^{j-1})] \Big\|^2 \Big| \mathcal{F}_{j} } - \Big \| \frac{1}{n} \sum_{i=1}^{n} [\nabla f_{i} (x^{j}) - \nabla f_{i}(x^{j-1})] \Big \|^2 \notag\\
&  =& \Exp{ \Big\|  \sum_{i \in I_{j}} \xi_i \Big\|^2 \Big| \mathcal{F}_{j} } - \Big \| \frac{1}{n} \sum_{i=1}^{n} \xi_i \Big \|^2\notag \\
&  \overset{\eqref{lem:upperv}}{\leq} & \sum_{i=1}^{n} v_ip_i \| \xi_i \|^2 \notag\\
&  = & \sum_{i=1}^{n}\frac{v_ip_i}{p_i^2n^2} \| \nabla f_{i} (x^{j}) - \nabla f_{i}(x^{j-1}) \|^2 \notag\\
& \overset{\eqref{L_i-smooth},\eqref{def_K}}{\leq} &\frac{1}{b} K \eta^2 \| v^{j-1} \|^2. \notag
\end{eqnarray}

Hence, by taking expectation, we have
\begin{eqnarray}
\Exp{\| v^j - v^{j-1} \|^2 } - \Exp{ \| \nabla f(x^{j}) - \nabla f(x^{j-1}) \|^2} \leq \frac{1}{b} K \eta^2 \Exp{[ \| v^{j-1} \|^2}. \notag
\end{eqnarray}

By Lemma \ref{lem:var_diff_mb}, for $t \geq 1$, 
\begin{eqnarray}
\Exp{\| \nabla f(x^{t}) - v^{t} \|^2 }
& = & \sum_{j = 1}^{t}\Exp{ \| v^{j} - v^{j-1} \|^2 } - \sum_{j = 1}^{t} \Exp{\| \nabla f(x^{j}) - \nabla f(x^{j-1}) \|^2 } \notag\\
 & \leq &\frac{1}{b}K \eta^2 \sum_{j = 1}^{t} \Exp{ \| v^{j-1} \|^2}.  \notag
\end{eqnarray}

This completes the proof. 
\end{proof}

\subsubsection*{Proof of Theorem \ref{thm:nonconvex_01_mb}}
\begin{proof}
By Lemma \ref{lem:var_diff_mb_02}, we have
\begin{eqnarray}
\Exp{ \| \nabla f(x^{t}) - v^{t} \|^2 } \leq \frac{1}{b} K\eta^2 \sum_{j=1}^{t} \Exp{\| v^{j-1} \|^2}.  \notag
\end{eqnarray}

Note that $\| \nabla f(x^{0}) - v^{0} \|^2 = 0$. Hence, by summing over $t = 0,\dots,m$ ($m \geq 1$), we have
\begin{eqnarray}
\sum_{t=0}^{m} \Exp{\| v^{t} - \nabla f(x^{t}) \|^2} &\leq & \frac{1}{b}K \eta^2 \Big[ m  \Exp{\|v^{0} \|^2} \notag \\
&& \quad+ (m-1)\Exp{\|v^{1} \|^2} + \dots +\Exp{\|v^{m-1}\|^2} \Big ].  
\end{eqnarray}

We have
\begin{eqnarray}
\label{eq:equal_zero_mb}
&& \sum_{t=0}^{m}\Exp{ \| \nabla f(x^{t}) - v^{t} \|^2 }  - ( 1 - \bar{L}\eta ) \sum_{t=0}^{m} \Exp{ \| v^{t} \|^2 } \notag\\
 && \leq  \frac{1}{b} K \eta^2 \Big[ m  \Exp{\|v^{0} \|^2} + (m-1) \Exp{\|v^{1} \|^2 }+ \dots + \Exp{\|v^{m-1}\|^2} \Big ]\notag\\ 
 && \quad- (1 - \bar{L}\eta) \Big [ \Exp{\|v^{0} \|^2} + \Exp{\|v^{1} \|^2} + \dots + \Exp{\|v^{m}\|^2 } \Big ] \notag\\
 && \leq  \Big[\frac{1}{b} K\eta^2 m - (1 - \bar{L}\eta) \Big] \sum_{t=1}^{m}\Exp{ \| v^{t-1} \|^2 } \overset{\eqref{eta_mb}}{\leq} 0
\end{eqnarray}

since 
\begin{eqnarray}
\eta = \frac{2}{\bar{L}\left(\sqrt{1 + \frac{4K m}{\bar{L}^2b}} + 1\right)}\notag
\end{eqnarray}

is a root of equation 
\begin{eqnarray}
\frac{1}{b}  K\eta^2 m - (1 - \bar{L}\eta) = 0. \notag
\end{eqnarray}

Therefore, by Lemma \ref{lem_main_derivation_mb}, we have
\begin{eqnarray}
\sum_{t=0}^{m}\Exp{ \| \nabla f(x^{t})\|^2 }& \leq &\frac{2}{\eta} [ f(x^{0}) - f(x^{*})] + \sum_{t=0}^{m} \Exp{ \| \nabla f(x^{t}) - v^{t} \|^2 }  \notag \\
&& \quad
 - ( 1 - \bar{L}\eta ) \sum_{t=0}^{m} \Exp{ \| v^{t} \|^2 }\notag\\
 & \overset{\eqref{eq:equal_zero_mb}}{\leq} \frac{2}{\eta}& [ f(x^{0}) - f(x^{*})].  \notag
\end{eqnarray}

If $x_a$ is chosen uniformly at random from $\{x^t\}_{t=0}^{m}$, then 
\begin{eqnarray}
\Exp{\| \nabla f(x_a)\|^2 } = \frac{1}{m+1}\sum_{t=0}^{m} \Exp{\| \nabla f(x^{t})\|^2 ] \leq \frac{2}{\eta(m+1)} [ f(x^{0}) - f(x^{*})}.  \notag
\end{eqnarray}
This concludes the proof.
\end{proof}

\section{One Sample Importance Sampling}

\subsection{SVRG}

\begin{algorithm}[hb]\small
   \caption{SVRG$\left(x^0,T, m, \{p_i\}_{i=0}^{n}, \eta \right)$}
   \label{alg:svrg_imp}
\begin{algorithmic}[1]
   \STATE {\bfseries Input:} $\tilde{x}^0 = x^0_m = x^0 \in \R^d$,  epoch length $m$, step sizes $\{\eta_i > 0\}_{i=0}^{m-1}$, $S = \lceil T/m \rceil$
   \FOR{$s=0$ {\bfseries to} $S-1$}
   \STATE $x^{s+1}_0 = x^{s}_m$
   \STATE $g^{s+1} = \frac{1}{n} \sum_{i=1}^n \nabla f_{i}(\tilde{x}^{s})$
   \FOR{$t=0$ {\bfseries to} $m-1$}
   \STATE With $\{p_i\}_{i=0}^{n} $ randomly pick $i_t$ from $\{1, \dots, n\}$ 
   \STATE $v_t^{s+1} = \frac{1}{np_{i_t}}( \nabla f_{i_t}(x^{s+1}_t) - \nabla f_{i_t}(\tilde{x}^{s})) + g^{s+1}$ 
   \STATE $x^{s+1}_{t+1} = x^{s+1}_{t} - \eta v_t^{s+1} $
   \ENDFOR
   \STATE $\tilde{x}^{s+1} = x_{m}^{s+1}$
   \ENDFOR
   \STATE {\bfseries Output:} Iterate $x_a$ chosen uniformly random from $\{\{x^{s+1}_t\}_{t=0}^{m}\}_{s=0}^{S}$.
\end{algorithmic}
\end{algorithm}

In this section, we introduce \texttt{SVRG} algorithm with batch size equal to $1$.

\begin{theorem}
   Let $c_m = 0$, $\eta = \eta > 0$, $\beta = \beta > 0$, and $c_{t} = c_{t+1}(1 + \eta\beta + 2K\eta^2 ) +  K\eta^2\bar{L}$ such that $\Gamma_t > 0$ for $0\leq t\leq m-1$. Define the quantity $\gamma_n \eqdef \min_t\Gamma_t$. 
  Further, let $T$ be a multiple of $m$. Then for the output $x_a$ of Algorithm~\ref{alg:svrg_imp} we have
  \begin{eqnarray}
    \Exp{\|\nabla f(x_a)\|^2} \leq \frac{f(x^{0}) - f(x^*)}{T\gamma_n},
  \end{eqnarray}
  where $x^*$ is an optimal solution to \eqref{eq:1} and $\Gamma_t = \bigl(\eta - \frac{c_{t+1}\eta}{\beta} - \eta^2\bar{L} - 2c_{t+1}\eta^2\bigr)$.
  \label{thm:nonconvex-inter_is}
\end{theorem}

\begin{theorem}
 Let $\eta = \bar{L}\mu_0/(Kn^{\frac23})$ ($0 < \mu_0 < 1$), $\beta = \bar{L}/n^{\frac13}$, $m = \lfloor K n/(3\bar{L}^2\mu_0) \rfloor$ and $T$ is some multiple of $m$. Then there exists universal constants $\mu_0, \nu > 0$ such that we have the following: $\gamma_n \geq \frac{\bar{L}\nu}{K}n^{\frac23}$ in Theorem~\ref{thm:nonconvex-inter_is} and
  \begin{eqnarray}
    \Exp{\|\nabla f(x_a)\|^2} &\leq \frac{Kn^{\frac23} [f(x^{0}) - f(x^*)]}{\bar{L}T\nu},
  \end{eqnarray} 
  where $x^*$ is an optimal solution to the problem in~\eqref{eq:1} and $x_a$ is the output of Algorithm~\ref{alg:svrg_imp}.
  \label{thm:nonconvex-gen_is}
\end{theorem}

Comparing Theorem~\ref{thm:nonconvex-inter_is} 
to the previous result in \cite{reddi2016stochastic}, we can see improvement in constant, if we assume different $L_i$-smooth constants for different functions. If the all $L_i$'s are the same then our result is the same as previous result for uniform sampling, because then $\alpha = \frac{n-1}{n-1}=1$.

\subsection{SAGA}

\begin{algorithm}[tb]\small
   \caption{SAGA$\left(x^0,T, \{p_i\}_{i=0}^{n},\eta\right)$}
   \label{alg:saga}
\begin{algorithmic}[1]
   \STATE {\bfseries Input:} $x^0 \in \R^d$, $\alpha_{i}^0 = x^0$ for $i \in [n]$,  number of iterations $T$, step size $\eta > 0$
   \STATE $g^{0} = \frac{1}{n} \sum_{i=1}^n \nabla f_{i}(\alpha_i^{0})$
   \FOR{$t=0$ {\bfseries to} $T-1$}
   \STATE  Randomly pick $i_t$ from $[n]$ with $\{p_i\}_{i=0}^{n}$
   \STATE  Randomly uniformly pick $i_t$ from $[n]$ 
   \STATE $v^t = \frac{1}{np_{i_t}} (\nabla f_{i_t}(x^t) - \nabla f_{i_t}(\alpha_{i_t}^{t})) + g^{t}$
   \STATE $x^{t+1} = x^{t} - \eta v^t$
   \STATE $\alpha_{j_t}^{t+1} = x^t$ and $\alpha_{j}^{t+1} = \alpha_{j}^{t}$ for $j \neq j_t$
   \STATE $g^{t+1} = g^t - \frac{1}{n}(\nabla f_{j_t}(\alpha_{j_t}^t) - \nabla f_{j_t}(\alpha_{j_t}^{t+1}))$
   \ENDFOR
   \STATE {\bfseries Output:} Iterate $x_a$ chosen uniformly random from $\{x^t\}_{t=0}^{T}$.
\end{algorithmic}
\end{algorithm}

Here, we provide similar analysis as for \texttt{SVRG} with the same result. We provide more generalized improved form of theorems which appeared in \cite{reddi2016fast}.

\begin{theorem}
\label{thm:nonconvex-inter}
Let $c_T = 0$, $\beta > 0$, and $c_{t} = c_{t+1}(1 - \tfrac{1}{n} + \eta\beta + 2K\eta^2) +  K\eta^2\bar{L}$ be such that $\Gamma_t > 0$ for $0 \leq t \leq T-1$. Define the quantity $\gamma_n \eqdef \min_{0 \leq t \leq T-1} \Gamma_t$. Then  the output $x_a$ of Algorithm~\ref{alg:saga} satisfies the bound
  \begin{eqnarray}
    \Exp{\|\nabla f(x_a)\|^2} \leq \frac{f(x^{0}) - f(x^*)}{T\gamma_n}, \notag 
  \end{eqnarray}
  where $x^*$ is an optimal solution to~\eqref{eq:1} and $\Gamma_{t} = \bigl(\eta - \tfrac{c_{t+1}\eta}{\beta} - \eta^2\bar{L} - 2c_{t+1}\eta^2\bigr)$.
\end{theorem}

\begin{theorem}
\label{thm:nonconvex-gen}
Let $\eta = \bar{L}/(3Kn^{2/3})$ and $\beta = \bar{L}/n^{1/3}$. Then,  $\gamma_n \geq \frac{\bar{L}}{12Kn^{2/3}}$ and we have the bound
\begin{eqnarray}
\Exp{\|\nabla f(x_a)\|^2} &\leq \frac{12Kn^{2/3} [f(x^{0}) - f(x^*)]}{\bar{L}T},\notag 
\end{eqnarray} 
where $x^*$ is an optimal solution to the problem in~\eqref{eq:1} and $x_a$ is the output of Algorithm~\ref{alg:saga}.
\end{theorem}

We can see that exactly same conclusions apply here as for \texttt{SVRG} and results can be interpreted in the same way.

\section{\texttt{SARAH:} Convex Case} \label{sec:SARAH-convex}

\subsection{Main result}

Consider Algorithm~\ref{alg:sarah-is}, which is an arbitrary sampling  variant of the \texttt{SARAH} method..

\begin{algorithm}
   \caption{\texttt{SARAH}}
      \label{alg:sarah-is}
\begin{algorithmic}[1]
   \STATE {\bfseries Parameters:} the learning rate $\eta > 0$ and the inner loop size $m$.
   \STATE {\bfseries Initialize:} $\tilde{x}_0$
   \STATE {\bfseries Iterate:}
   \FOR{$s=1,2,\dots$}
   \STATE $x_0 = \tilde{x}_{s-1}$
   \STATE $v^0 = \frac{1}{n}\sum_{i=1}^{n} \nabla f_i(x^0)$
   \STATE $x_1 = x_0 - \eta v^0$
   \STATE {\bfseries Iterate:}
   \FOR{$t=1,\dots,m-1$}
   \STATE Sample $i_{t}$ at random from $[n]$ with probability $\{ p_i\}_{i=1}^n$
   \STATE $v^{t} = \frac{1}{np_i}(\nabla f_{i_{t}} (x^{t}) - \nabla f_{i_{t}}(x^{t-1})) + v^{t-1}$
   \STATE $x_{t+1} = x^{t} - \eta v^{t}$
   \ENDFOR
   \STATE Set $\tilde{x}_s = x^{t}$ with $t$ chosen uniformly at random from $\{0,1,\dots,m\}$
   \ENDFOR
\end{algorithmic}
\end{algorithm}

Note, that only $10$-th and $11$-th row are changed comparing to classic \texttt{SARAH}  algorithm presented in \cite{SARAH}. We do not sample uniformly anymore and also in the $11$-th row of Algorithm~\ref{alg:sarah-is}, where we use factor $\frac{1}{np_i}$ in order to stay unbiased in outer cycle. 

Then using similar analysis used in \cite{SARAH} and additional lemmas we can prove following theorems with $p_i$ in Algorithm~\ref{alg:sarah-is} to be $\frac{L_i}{\sum_{j=1}^nL_i}$

\begin{theorem}\label{lem_bouned_moment_stronglyconvexP}
Suppose that $f_i(x)$ are $L_i$-smooth and convex, $f(x)$ is $\mu$ strongly convex. Consider $v^{t}$ defined in  \texttt{SARAH} (Algorithm~\ref{alg:sarah-is}) with $\eta < 2/\bar{L}$, where $\bar{L} = \frac1n \sum_{j=1}^nL_i$. Then, for any $t\geq 1$,
\begin{eqnarray}
\Exp{\|v^{t}\|^2}
 &\leq & \left[ 1 - \left(\tfrac{2}{\eta \bar{L}} - 1 \right) \mu^2 \eta^2  \right] \Exp{\|v^{t-1}\|^2}
 \notag\\
 &\leq &\left[ 1 - \left(\tfrac{2}{\eta \bar{L}} - 1 \right) \mu^2 \eta^2  \right]^{t} \Exp{ \| \nabla f(x^{0}) \|^2 }.\notag
\end{eqnarray}
\end{theorem}
By choosing $\eta=\Ocal(1/\bar{L})$, we obtain the linear  convergence of $\|v^t\|^2$ in expectation with the rate $(1-1/\kappa^2)$, where $\kappa	= \frac{\bar{L}}{\mu}$ is condition number, This is improvement over previous result in \cite{SARAH}, because of $\frac{\bar{L}}{\mu}\leq \frac{L_{\max}}{\mu}$. Below we show that a better convergence rate could be obtained under a  stronger convexity assumption for each single $f_i(x)$. 

\begin{theorem}\label{thm:bound_moment}
	Suppose that $f_i(x)$ are $L_i$-smooth and $\mu$ strongly convex.  Consider $v^{t}$ defined by  in \texttt{SARAH} (Algorithm \ref{alg:sarah-is}) with $\eta \leq 2/(\mu+\bar{L})$. Then the following bound holds, $\forall\ t\geq 1$, 
\begin{eqnarray}
\Exp{\| v^{t} \|^2 }
& \leq & \left( 1 - \tfrac{2 \mu \bar{L} \eta}{\mu +\bar{L}} \right) \Exp{  \|v^{t-1} \|^2 }
\notag\\
& \leq &\left(1 - \tfrac{2\mu \bar{L} \eta}{\mu + \bar{L}} \right)^{t} \Exp{ \| \nabla f(x^{0}) \|^2 }. \notag
\end{eqnarray}
\end{theorem}

By setting $\eta=\Ocal(1/\bar{L})$, we derive the linear convergence with the rate of $(1-1/\kappa)$, where $\hat{\kappa} = \frac{\bar{L}}{\mu}$  which is an improvement over the previous result of \cite{SARAH}, because if we take the optimal stepsize $\nu = \tfrac{2}{\mu + \bar{L}}$ than we can easily prove that $\tfrac{2\mu \bar{L} \eta}{\mu + \bar{L}}$ is greater than $\tfrac{2\mu L_{\max} \eta}{\mu + L_{\max}}$, with optimal step size, where $L_{\max} = \max_i\{L_i\}$.

\subsection{Lemmas}

We start with modification of lemmas in \cite{SARAH}, which we later use in the proofs of Theorem~\ref{thm:bound_moment} and Theorem~\ref{lem_bouned_moment_stronglyconvexP}. 
The first Lemma 
\ref{lem_main_derivation}
bounds the  sum of expected values of
$\|\nabla f(x^t)\|^2$.
The second, Lemma \ref{lem:var_diff_01},
bounds $\Exp{ \| \nabla f(x^{t}) - v^{t} \|^2 }$. 

\begin{lemma}\label{lem_main_derivation}
Suppose that $f_i(x)$'s are $L_i$-smooth. Consider \texttt{SARAH} (Algorithm \ref{alg:sarah-is}). Then, we have 
\begin{eqnarray}
& &\sum_{t=0}^{m} \Exp{ \| \nabla f(x^{t})\|^2 }  \leq \frac{2}{\eta} \Exp{ f(x^{0}) - f(x^{*})}
\label{eq:001_1} 
\notag\\&&\quad+ \sum_{t=0}^{m} \Exp{\| \nabla f(x^{t}) - v^{t} \|^2 }  
 - ( 1 - \bar{L}\eta ) \sum_{t=0}^{m} \Exp{ \| v^{t} \|^2 }.
\end{eqnarray}
\end{lemma}
\begin{lemma}\label{lem:var_diff_01}
Suppose that $f_i(x)$'s are $L_i$-smooth. Consider  \texttt{SARAH}  (Algorithm \ref{alg:sarah-is}). Then for any $t\geq 1$, 
\begin{eqnarray}
\Exp{ \| \nabla f(x^{t}) - v^{t} \|^2 } 
&=& \sum_{j = 1}^{t} \Exp{ \| v^{j} - v^{j-1} \|^2 }  
- \sum_{j = 1}^{t} \Exp{ \| \nabla f(x^{j}) - \nabla f(x^{j-1}) \|^2 }. \notag
\end{eqnarray}
\end{lemma}

\begin{lemma}\label{lem_bound_var_diff_str_02}
Suppose that $f_i(x)$'s  are $L_i$-smooth and convex. Consider  \texttt{SARAH} (Algorithm \ref{alg:sarah-is}) with $\eta < 2/\bar{L}$. Then we have that for any $t\geq 1$, 
\begin{eqnarray}
\Exp{ \| \nabla f(x^{t}) - v^{t} \|^2 } 
&\leq & \frac{\eta \bar{L}}{2 - \eta \bar{L}} \Big[ \Exp{ \|v^{0} \|^2} - \Exp{\| v^{t} \|^2} \Big]\notag
\\
&\leq&  \frac{\eta \bar{L}}{2 - \eta \bar{L}}  \Exp{\|v^{0} \|^2},\notag
\label{eq:bound1}
\end{eqnarray}
where $\bar{L} = \frac{1}{n}\sum_{i=1}^n L_i$.
\end{lemma}

\subsubsection*{Proof of Lemma \ref{lem_main_derivation}}
\begin{proof}
By Lemma~\ref{lem:liptschitz} and $x^{t+1} = x^{t} - \eta v^{t}$, we have
\begin{eqnarray}
\Exp{f(x^{t+1})} & \leq & \Exp{ f(x^{t})} - \eta \Exp{\nabla f(x^{t})^\top v^{t}} 
+ \frac{\bar{L}\eta^2}{2} \Exp{ \| v^{t} \|^2 } \notag
\\
& = &\Exp{ f(x^{t})}- \frac{\eta}{2} \Exp{\| \nabla f(x^{t})\|^2 }
+ \frac{\eta}{2} \Exp{ \| \nabla f(x^{t}) - v^{t} \|^2 } \notag \\
&& \quad
- \left( \frac{\eta}{2} - \frac{\bar{L}\eta^2}{2} \right) \Exp{ \| v^{t} \|^2 },\notag
\end{eqnarray}
where the last equality follows from the fact
$a^\top b = \frac{1}{2}\left[\|a\|^2 + \|b\|^2 - \|a-b\|^2\right].$

By summing over $t = 0,\dots,m$, we have
\begin{eqnarray}
\Exp{ f(x_{m + 1})} & \leq & \Exp{ f(x^{0})} - \frac{\eta}{2} \sum_{t=0}^{m} \Exp{ \| \nabla f(x^{t})\|^2 } + \frac{\eta}{2} \sum_{t=0}^{m} \Exp{ \| \nabla f(x^{t}) - v^{t} \|^2 }  \notag \\
&& \quad
- \left( \frac{\eta}{2} - \frac{\bar{L}\eta^2}{2} \right) \sum_{t=0}^{m} \Exp{ \| v^{t} \|^2 },  \notag
\end{eqnarray}
which is equivalent to ($\eta>0$):
\begin{eqnarray}
\sum_{t=0}^{m} \Exp{ \| \nabla f(x^{t})\|^2 }  & \leq &\frac{2}{\eta} \Exp{ f(x^{0}) - f(x_{m+1})}+\notag \\
&&\quad \sum_{t=0}^{m} \Exp{ \| \nabla f(x^{t}) - v^{t} \|^2 }  
 - ( 1 - \bar{L}\eta ) \sum_{t=0}^{m} \Exp{ \| v^{t} \|^2 } \notag\\
& \leq &\frac{2}{\eta} \Exp{ f(x^{0}) - f(x^{*})} + \sum_{t=0}^{m} \Exp{ \| \nabla f(x^{t}) - v^{t} \|^2 } \notag \\ && \quad
 - ( 1 - \bar{L}\eta ) \sum_{t=0}^{m} \Exp{\| v^{t} \|^2 },    \notag
\end{eqnarray}
where the last inequality follows since $x^{*}$ is a global minimizer of \eqref{eq:1}. 
\end{proof}

\subsubsection*{Proof of Lemma \ref{lem:var_diff_01}}
\begin{proof}
Let $\mathcal{F}_{j}$ be $\sigma$ algebra that contains all the information of $x^{0},\dots,x^{j}$ as well as $v^0,\dots,v^{j-1}$. For $j \geq 1$, we have
\begin{eqnarray}
&&\Exp{ \| \nabla f(x^{j}) - v^{j} \|^2 | \mathcal{F}_{j} } 
 =  \notag \\
&& \Exp{ \| [\nabla f(x^{j-1}) - v^{j-1} ] + [ \nabla f(x^{j}) - \nabla f(x^{j-1}) ] - [ v^{j} - v^{j-1} ] \|^2 | \mathcal{F}_{j} } \notag
\\
&& = \| \nabla f(x^{j-1}) - v^{j-1} \|^2 + \| \nabla f(x^{j}) - \nabla f(x^{j-1}) \|^2 \notag\\
&&\quad+ \Exp{ \| v^{j} - v^{j-1}  \|^2 | \mathcal{F}_{j} }\notag
\\
&&\quad + 2 ( \nabla f(x^{j-1}) - v^{j-1} )^\top ( \nabla f(x^{j}) - \nabla f(x^{j-1}) )\notag \\
&&\quad - 2 ( \nabla f(x^{j-1}) - v^{j-1} )^\top \Exp{v^{j} - v^{j-1} | \mathcal{F}_{j} } \notag \\
&&\quad - 2 ( \nabla f(x^{j}) - \nabla f(x^{j-1}) )^\top \Exp{v^{j} - v^{j-1} | \mathcal{F}_{j} } \notag
\\
&& = \| \nabla f(x^{j-1}) - v^{j-1} \|^2 - \| \nabla f(x^{j}) - \nabla f(x^{j-1}) \|^2  \notag \\
&& \quad +\Exp{ \| v^{j} - v^{j-1}  \|^2 | \mathcal{F}_{j} }, \notag
\end{eqnarray}
where the last equality follows from
\begin{eqnarray}
& \Exp{v^{j} - v^{j-1} | \mathcal{F}_{j} }= \Exp{\frac{1}{np_{i_j}}\left( \nabla f_{i_{j}}(x^{j}) - \nabla f_{i_{j}}(x^{j-1})\right) | \mathcal{F}_{j} } 
= \nabla f(x^{j}) - \nabla f(x^{j-1}). \notag 
\end{eqnarray}

By taking expectation for the above equation, we have
\begin{eqnarray}
\Exp{ \| \nabla f(x^{j}) - v^{j} \|^2 } &=& \Exp{\| \nabla f(x^{j-1}) - v^{j-1} \|^2 } - \Exp{ \| \nabla f(x^{j}) - \nabla f(x^{j-1}) \|^2 }\notag \\
&& \quad + \Exp{\| v^{j} - v^{j-1} \|^2 }. \notag 
\end{eqnarray}

Note that $\| \nabla f(x^{0}) - v^{0} \|^2 = 0$. By summing over $j = 1,\dots,t\ (t\geq 1)$, we have
\begin{eqnarray}
& \Exp{ \| \nabla f(x^{t}) - v^{t} \|^2 }  = \sum_{j = 1}^{t} \Exp{ \| v^{j} - v^{j-1} \|^2 }  - \sum_{j = 1}^{t} \Exp{ \| \nabla f(x^{j}) - \nabla f(x^{j-1}) \|^2 }. \notag 
\end{eqnarray}
\end{proof}
\subsubsection*{Proof of Lemma \ref{lem_bound_var_diff_str_02}}
\begin{proof}
For $j \geq 1$, we have
\begin{eqnarray}
 \Exp{\| v^{j} \|^2 | \mathcal{F}_{j} } &=& \quad \Exp{\| v^{j-1} - \frac{1}{np_{i_j}}(\nabla f_{i_{j}}(x^{j-1}) - \nabla f_{i_{j}}(x^{j}) ) \|^2 | \mathcal{F}_{j} } \notag 
\\ 
&=& \|v^{j-1} \|^2 + \Exp{\frac{1}{n^2p_{i_j}^2} \| \nabla f_{i_{j}}(x^{j-1}) - \nabla f_{i_{j}}(x^{j}) \|^2 | \mathcal{F}_{j}}\notag \\
&& \quad  - \Exp{\tfrac{2}{\eta np_{i_j}}(\nabla f_{i_{j}}(x^{j-1}) - \nabla f_{i_{j}}(x^{j}))^\top (x^{j-1} - x^{j}) | \mathcal{F}_{j}} \notag 
 \\ 
& \overset{\eqref{L_i-smooth}}{\leq} &\|v^{j-1} \|^2  + \Exp{\frac{1}{n^2p_{i_j}^2} \| \nabla f_{i_{j}}(x^{j-1}) - \nabla f_{i_{j}}(x^{j}) \|^2 | \mathcal{F}_{j}} \notag \\
&& \quad - \Exp{\tfrac{2}{L_{i_j} \eta np_{i_j}} \| \nabla f_{i_{j}}(x^{j-1}) - \nabla f_{i_{j}}(x^{j}) \|^2 | \mathcal{F}_{j}} 
\notag 
\\ 
& = &\|v^{j-1} \|^2 + \left(1 - \tfrac{2}{\eta \bar{L}}\right) \Exp{ \left\|\frac{1}{np_{i_j}}\left( \nabla f_{i_{j}}(x^{j-1}) - \nabla f_{i_{j}}(x^{j})\right) \right\|^2 | \mathcal{F}_{j} } \notag \\
& =& \|v^{j-1} \|^2 + \left(1 - \tfrac{2}{\eta \bar{L}}\right) \Exp{ \| v^{j} - v^{j-1} \|^2 | \mathcal{F}_{j} },\notag 
\end{eqnarray}
 The consequent equality follows from definition of $p_i$'s and the last equality follows from definition of  \texttt{SARAH} .
Taking expectation, we get 
\begin{eqnarray}
\Exp{\| v^{j} - v^{j-1} \|^2}
\leq \frac{\eta \bar{L}}{2 - \eta \bar{L}} \Big[ \Exp{ \|v^{j-1} \|^2} - \Exp{\| v^{j} \|^2} \Big],\notag 
\end{eqnarray}
when $\eta < 2/{\bar{L}}$.

\quad By summing the above inequality over $j = 1,\dots, t\ (t\geq 1)$, we have
\begin{equation}\label{eq:sumover}
\sum_{j=1}^{t}\Exp{\| v^{j} - v^{j-1} \|^2}
\leq \frac{\eta \bar{L}}{2 - \eta \bar{L}} \Big[ \Exp{ \|v^{0} \|^2} - \Exp{\| v^{t} \|^2} \Big].  
\end{equation}

By Lemma \ref{lem:var_diff_01}, we have
\begin{eqnarray}
 \Exp{ \| \nabla f(x^{t}) - v^{t} \|^2 } & \leq \sum_{j = 1}^{t} \Exp{ \| v^{j} - v^{j-1} \|^2 }  \overset{\eqref{eq:sumover}}{\leq}  \frac{\eta \bar{L}}{2 - \eta \bar{L}} \Big[ \Exp{ \|v^{0} \|^2} - \Exp{\| v^{t} \|^2} \Big]. \notag 
\end{eqnarray}
\end{proof}

\subsubsection*{Proof of Theorem~\ref{lem_bouned_moment_stronglyconvexP}}
\begin{proof}
For $t \geq 1$, we have
\begin{eqnarray}
\| \nabla f(x^{t}) - \nabla f(x^{t-1})\|^2 &=& \Big\| \frac{1}{n} \sum_{i=1}^{n} [ \nabla f_i(x^{t}) - \nabla f_i(x^{t-1})  ] \Big\|^2 \notag \\
&=& \Big\| \sum_{i=1}^{n} p_i \frac{1}{np_i}[ \nabla f_i(x^{t}) - \nabla f_i(x^{t-1})  ] \Big\|^2\notag  \\
& \overset{\eqref{Jensen}}{\leq} &\sum_{i=1}^{n} p_i  \| \tfrac{1}{np_i}\left(\nabla f_i(x^{t}) - \nabla f_i(x^{t-1}) \right) \|^2 \notag \\
& = &\Exp{ \left\|\tfrac{1}{np_i} \left( \nabla f_{i_{t}}(x^{t}) - \nabla f_{i_{t}}(x^{t-1})\right) \right\|^2 | \mathcal{F}_{t} }.  \label{eq:prthm001}
\end{eqnarray}

Using the proof of Lemma \ref{lem_bound_var_diff_str_02}, for $t \geq 1$, we have
\begin{eqnarray}
 \Exp{\| v^{t} \|^2 | \mathcal{F}_{t}} 
& \leq & \|v^{t-1} \|^2 + \left(1 - \tfrac{2}{\eta \bar{L}}\right) \Exp{\| \nabla f_{i_{t}}(x^{t-1}) - \nabla f_{i_{t}}(x^{t}) \|^2 | \mathcal{F}_{t}} \notag \\
& \overset{\eqref{eq:prthm001}}{\leq } &\|v^{t-1} \|^2 + \left(1 - \tfrac{2}{\eta \bar{L}}\right) \| \nabla f(x^{t}) - \nabla f(x^{t-1})\|^2 \notag \\
& \leq &\|v^{t-1} \|^2 + \left(1 - \tfrac{2}{\eta \bar{L}}\right) \mu^2 \eta^2 \|v^{t-1}\|^2.\notag 
\end{eqnarray}

Note that $1 - \tfrac{2}{\eta \bar{L}} < 0$ since $\eta < 2/\bar{L}$. The last inequality follows by the strong convexity of $f$, that is, $\mu \|x^{t} - x^{t-1}\| \leq \| \nabla f(x^{t}) - \nabla f(x^{t-1})\|$ and the fact that $x^{t} = x^{t-1} - \eta v^{t-1}$. By taking the expectation and applying recursively, we have
\begin{eqnarray}
\Exp{\| v^{t} \|^2} & \leq& \left[ 1 - \left(\tfrac{2}{\eta \bar{L}} - 1 \right) \mu^2 \eta^2  \right] \Exp{\| v^{t-1} \|^2} \notag \\
& \leq& \left[ 1 - \left(\tfrac{2}{\eta \bar{L}} - 1 \right) \mu^2 \eta^2  \right]^{t} \Exp{\| v^{0} \|^2}  \notag \\
& = &\left[ 1 - \left(\tfrac{2}{\eta \bar{L}} - 1 \right) \mu^2 \eta^2  \right]^{t} \Exp{\| \nabla f(x^{0}) \|^2}.  \notag
\end{eqnarray}
\end{proof}

\subsubsection*{Proof of Theorem~\ref{thm:bound_moment}}
\begin{proof}
We obviously have $\Exp{\|v^0\|^2 | \mathcal{F}_{0}} = \|\nabla f(x_0)\|^2$. For $t \geq 1$, we have
\begin{eqnarray}
\Exp{\| v^{t} \|^2 | \mathcal{F}_{t}} 
& = &\Exp{\| v^{t-1} - \tfrac{1}{np_{i_t}}(\nabla f_{i_{t}}(x^{t-1}) - \nabla f_{i_{t}}(x^{t}) ) \|^2 | \mathcal{F}_{t}} \notag \\
& = & \|v^{t-1} \|^2  + \Exp{\tfrac{1}{n^2p_{i_t}^2} \| \nabla f_{i_{t}}(x^{t-1}) - \nabla f_{i_{t}}(x^{t}) \|^2\mathcal{F}_{t} } \notag  \\  
&&  - \Exp{\tfrac{2}{\eta np_{i_t}}(\nabla f_{i_{t}}(x^{t-1})- \nabla f_{i_{t}}(x^{t}))^\top (x^{t-1} - x^{t}) | \mathcal{F}_{t} }\notag  \\
& = &\|v^{t-1} \|^2    + \Exp{ \| \tfrac{1}{np_{i_t}}\left(\nabla f_{i_{t}}(x^{t-1}) - \nabla f_{i_{t}}(x^{t})\right) \|^2| \mathcal{F}_{t} } \notag \\
&&   - \tfrac{2}{\eta}(\nabla f(x^{t-1}) - \nabla f(x^{t}))^\top (x^{t-1} - x^{t}) \notag \\
&\overset{\eqref{L_bar-smooth},\eqref{wrweq}}{\leq} &  \|v^{t-1} \|^2    +\Exp{\| \tfrac{1}{np_{i_t}}\left(\nabla f_{i_{t}}(x^{t-1}) - \nabla f_{i_{t}}(x^{t})\right) \|^2| \mathcal{F}_{t} }   \notag \\
&&\quad - \tfrac{2\mu\bar{L}\eta}{\mu + \bar{L}}\|v^{t-1} \|^2 -  \tfrac{2}{\eta(\mu + \bar{L})}\| \nabla f_{i_{t}}(x^{t-1}) - \nabla f_{i_{t}}(x^{t}) \|^2  \notag \\
& \leq &
\left( 1 - \tfrac{2 \mu \bar{L} \eta}{\mu + \bar{L}} \right)  \|v^{t-1} \|^2   \notag \\
&& 
+  \Exp{ \| \tfrac{1}{np_{i_t}}\left(\nabla f_{i_{t}}(x^{t-1}) - \nabla f_{i_{t}}(x^{t})\right) \|^2| \mathcal{F}_{t} }- \| \nabla f(x^{t-1}) - \nabla f(x^{t}) \|^2 \notag \\
&=& \left( 1 - \tfrac{2 \mu \bar{L} \eta}{\mu + \bar{L}} \right)  \|v^{t-1} \|^2   \notag \\
&& -\Exp{\| \tfrac{1}{np_{i_t}}\left(\nabla f_{i_{t}}(x^{t-1}) - \nabla f_{i_{t}}(x^{t})\right) -  \nabla f(x^{t-1}) - \nabla f(x^{t})\|^2|\mathcal{F}_{t} }\notag  \\
& \leq &  \left( 1 - \tfrac{2 \mu \bar{L} \eta}{\mu + \bar{L}} \right)  \|v^{t-1} \|^2,  
\label{eqasfewafaw}
\end{eqnarray}
where in the first two equalities, we used definition of  \texttt{SARAH} . The first inequality follows from fact that $f(x)$ is $\bar{L}$-smooth and $\mu$ strongly convex, thus following inequality holds (inequality from \cite{nesterov2013introductory})
\begin{equation}
\label{wrweq}
(\nabla f(x) - \nabla f(x'))^\top (x - x') \geq \frac{\mu \bar{L}}{\mu + \bar{L}} \|x - x'\|^2  
 + \frac{1}{\mu + \bar{L}}\| \nabla f(x) - \nabla f(x') \|^2, 
\end{equation}

 The second one uses assumption that $\eta \leq \tfrac{2}{\mu + \bar{L}}$, thus  $\eta =  \tfrac{2}{\mu + \bar{L}}$ is  optimal step size under this analysis. By taking the expectation and applying recursively, the desired result is achieved. 
\end{proof}

\section{Technical Lemmas}

\begin{lemma}
\label{lem:liptschitz}
  Let $f_{i}$'s be function, which are $L_i$-smooth, then $f(x) = \frac1n \sum_{i=1}^n f_i(x)$ is $\bar{L}$-smooth, where $\bar{L} = \frac1n \sum_{i=1}^n L_i$.
\end{lemma}

\begin{proof}
For each function $f_i$ we have by definition of $L_i$-smoothness, $\forall x,y \in \R^d$
\begin{equation}
\label{L_i-smooth}
f_i(x) \leq f_i(y) + \nabla f_i(y)^\top (x-y) + \frac{L_i}{2}\norm{x-y}^2 
\end{equation}
Summing through all $i$'s and dividing by $n$, we get
\begin{equation}
\label{L_bar-smooth}
f(x) \leq f(y) + \nabla f(y)^\top (x-y) + \frac{\bar{L}}{2}\norm{x-y}^2 
\end{equation}
\end{proof}

\begin{lemma}[Cauchy-Schwarz inequality] For all $x,y \in \R^d$ we have
\begin{equation}
\label{CS}
 |\langle x,y \rangle |\leq \|x \|\|y \|.
\end{equation}
\end{lemma}

\begin{lemma}[Young's inequality] For $a,b \in \R$ and $\beta >0$ we have
\begin{equation}
\label{Young}
 ab \leq \frac{a^2\beta}{2} + \frac{b^2}{2\beta}.
\end{equation}
\end{lemma}

\begin{lemma}[Jensen's inequality] Let $X$ be a random variable and $g(x)$ be a convex function. Then 
\begin{equation}
\label{Jensen}
 g(\Exp{X}) \leq \Exp{g(X)}.
\end{equation}
\end{lemma}

\end{document}